\documentclass[12pt,reqno]{amsart}
\usepackage{amsmath}
\usepackage{amssymb}
\usepackage{amsthm}
\usepackage{mathrsfs}
\usepackage{latexsym}
\usepackage{xcolor}
\usepackage{comment}

\newtheorem{theorem}{Theorem}[section]
\newtheorem{lemma}{Lemma}[section]
\newtheorem{corollary}{Corollary}[section]
\newtheorem{remark}{Remark}[section]
\newtheorem{definition}{Definition}[section]
\newtheorem{proposition}{Proposition}[section]
\newtheorem{example}{Example}[section]
\newtheorem{assumption}{Assumption}[section]
\numberwithin{equation}{section}
\newcommand{\bth}{\begin{theorem}}
\newcommand{\ethe}{\end{theorem}}

\newcommand{\bre}{\begin{remark}}
\newcommand{\ere}{\end{remark}}

\newcommand{\ble}{\begin{lemma}}
\newcommand{\ele}{\end{lemma}}
\newcommand{\bde}{\begin{definition}}
\newcommand{\ede}{\end{definition}}

\newcommand{\bco}{\begin{corollary}}
\newcommand{\eco}{\end{corollary}}

\newcommand{\bpr}{\begin{proposition}}
\newcommand{\epr}{\end{proposition}}

\newcommand{\bexer}{\begin{exercise}}
\newcommand{\eexer}{\end{exercise}}
\newcommand{\breh}{\begin{hint}}
\newcommand{\ereh}{\end{hint}}

\newcommand{\halmos}{\hfill \qed}

\newcommand{\bexam}{\begin{example}}
\newcommand{\eexam}{\end{example}}

\newcommand{\pr} {{\bf Proof.}}

\newcommand{\bfi}{\begin{fig}}
\newcommand{\efi}{\end{fig}}

\newcommand{\beao}{\begin{eqnarray*}}
\newcommand{\eeao}{\end{eqnarray*}\noindent}

\newcommand{\beam}{\begin{eqnarray}}
\newcommand{\eeam}{\end{eqnarray}\noindent}
\newcommand{\E}{\mathbf{E}}
\newcommand{\PP}{\mathbf{P}}
\newcommand{\nto}{n\to\infty}

\newcommand{\xto}{x\to\infty}

\newcommand{\bF}{\overline{F}}

\newcommand{\bG}{\overline{G}}
\newcommand{\bV}{\overline{V}}

\newcommand{\bbr}{{\mathbb R}}

\newcommand{\bbb}{{\mathbb B}}

\newcommand{\bbn}{{\mathbb N}}

\newcommand{\vep}{\varepsilon}

\allowdisplaybreaks[1]

\begin{document}
\title[Asymptotics for a nonstandard risk model]{Asymptotics for a nonstandard risk model with multivariate subexponential claims and constant interest force}

\author[ D.G. Konstantinides, C.D. Passalidis, H. Xu]{ Dimitrios G. Konstantinides, Charalampos  D. Passalidis, Hui Xu}

\address{Dept. of Statistics and Actuarial-Financial Mathematics,
University of the Aegean, Karlovassi, GR-83 200 Samos, Greece}
\address{School of Statistics and Data Science, Shanghai University of Finance and Economics, Shanghai, China}
\email{konstant@aegean.gr, \;sasd24009@sas.aegean.gr, \;xuhui@mail.shufe.edu.cn.}

\date{{\small \today}}

\begin{abstract}
In this paper, the asymptotic behavior of the entrance probability of discounted aggregate claims of a certain family of rare sets is studied, considering the finite and infinite time horizons. This multivariate risk model, driven by a common counting process, has a constant interest rate and allows the interdependence of claim vectors.  For the finite time horizon, the multivariate subexponential distribution of the common claim vector and the weak dependence structure of regression dependence are used. For the infinite time horizon, the claim vector is taken from a smaller distribution class, and the weak dependence structure is more general. Both results are derived under some additional assumptions on the moments of the counting process, which is fulfilled by all inhomogeneous renewal processes and many quasi-renewal processes, respectively. Moreover, the results are specialized to the multivariate regularly varying case, where more explicit results on the asymptotic behavior of the entrance probability of the discounted aggregate claims are derived. At the end of the paper, the results obtained are used to study the finite and infinite time horizon ruin problems of a risk model with eventual Brownian perturbations.
\end{abstract}

\maketitle
\vspace{3mm}



\section{Introduction} \label{sec.KPX.1}

\subsection{Model description}

In this paper, we consider an insurer operating $d$-lines of business, with $d \in \bbn$, each sharing a common counting process for claims. Let the sequence of random vectors $\{{\bf X}^{(i)},\,i \in \bbn\}$ represent claim vectors, where each ${\bf X}^{(i)}=\left(X_{1}^{(i)},\,\ldots,\, X_{d}^{(i)}\right)^{\top}$ and ${\bf X}^{(i)}$ follows a distribution defined on $\bbr_+^d := [0,\,\infty)^d$, the components of each vector may include zeros, but not all components can be zero.  Claims arrive at times $\{\tau_{i},\,i \in \bbn\}$, with the convention $\tau_0=0$, defining a counting process
\beam \label{eq.KPX.1.1}
N(t) := \sup \{i \in \bbn\;:\;\tau_i \leq t\}\,,
\eeam
for any $t \geq 0$, with $\sup \emptyset = 0$ by convention. The counting process $\{N(t)\,,\;t\geq 0\}$ is assumed to be a
general counting process, with a finite mean function
\beao
m(t) = \E[N(t)] = \sum_{i=1}^{\infty} \PP(\tau_i \leq t) < \infty\,,
\eeao
for any fixed $t\geq 0$. We restrict our analysis to the time interval $\Lambda :=\{t>0\;:\;m(t)> 0\}$, since the probability $P({\bf D}(T) \in x\,A)$ for $T\notin \Lambda$ are zero.
The insurer invests their surplus in risk free investments with a constant interest rate $r \geq 0$. Therefore, the insurer's discounted aggregated claims up to time $T \in \Lambda \setminus \{\infty\}$ are given by
\beam \label{eq.KPX.1.2}
{\bf D}(T) =\sum_{i=1}^{N(T)} {\bf X}^{(i)}\,e^{-r\,\tau_i}=\left(\sum_{i=1}^{N(T)} X_{1}^{(i)}\,e^{-r\,\tau_{i} },
\cdots,\sum_{i=1}^{N(T)} X_{d}^{(i)}\,e^{-r\,\tau_{i} }
\right)^{\top}\,,
\eeam
For the infinite time horizon, $T = \infty$, the discounted aggregate claims are defined similarly, with $N(\infty) = \infty$.

Our goal is to derive asymptotic estimations for the entrance probabilities of ${\bf D}(T)$ and ${\bf D}(\infty)$ into some rare set $x\,A$, where $A$ belongs to
a broad family of sets $\mathscr{R}$ as $x\to\infty$, i.e.,
\beam \label{eq.KPX.1.3}
\PP[{\bf D}(T) \in x\,A]\,, \qquad \PP[{\bf D}(\infty) \in x\,A]\,, \quad x\to\infty.
\eeam
These probabilities are related to ruin probabilities (for corresponding ruin sets), but they can also be used in broader applications of risk management and insurance solvency. For example, in the one-dimensional subcase with $A =(1,\,\infty)$, these probabilities provide tails for discounted aggregate claims in insurance portfolios, over both finite and infinite time horizons. We review the relevant literature related to these results in the next subsection. Notably, even in the one-dimensional case where $A=(1,\,\infty)$, the main results of this paper remain novel.

Let us present the following assumption, which will hold hereafter. We remind that the counting process $\{N(t)\,,\;t\geq 0\}$ from \eqref{eq.KPX.1.1} is a general counting process with mean function $m(t)$.

\begin{assumption} \label{ass.KPX.1.1}
The claim vectors $\{{\bf X}^{(i)}\,,\;i \in \bbn\}$ are identically distributed copies of the general random vector ${\bf X}$, following a distribution $F$ with support on  $\bbr_+^d$. We also assume that the $\{{\bf X}^{(i)},\,i \in \bbn\}$ and $\{\tau_{i},\,i \in \bbn\}$ are independent sequences.
\end{assumption}

\subsection{Overview of existing results} \label{subsection.1.2}

Our results aim to address two key objectives, which are of particular importance in actuarial practice:
\begin{enumerate}
\item
Non-standard driven multivariate risk models: This refers to moving beyond standard renewal processes, or quasi-renewal processes, in risk modeling.
\item
Interdependent claim vectors: We consider the possibility of both dependence within claim vectors and dependence across different claim vectors over time.
\end{enumerate}

In terms of non-standard counting processes, we look at models that don't necessarily follow the conventional renewal or quasi-renewal patterns. While many studies have tried to push the boundaries of traditional models in the context of one-dimensional risk models, such as those discussed in \cite{albrecher:asmussen:2006}, \cite{asmussen:schmidli:schmidt:1999}, \cite{bernackaite:siaulys:2017}, \cite{palmowski:pojer:thonhauser:2025}, these works primarily explore asymptotics for ruin probabilities in models without considering investments. Our work differs in that we build on \cite{cheng:xu:2020}, which considered a one-dimensional, inhomogeneous Poisson-driven risk model with risk free investments $(r>0)$ and derived asymptotic estimations for ruin probabilities over finite time.

In the more recent literature, attention has shifted to multivariate risk models, which are more practical for real-world applications, as seen in  \cite{chen:konstantinides:passalidis:2025}, \cite{chen:wang:wang:2013}, \cite{cheng:konstantinides:wang:2024}, \cite{gao:yang:2014}, \cite{jiang:wang:chen:xu:2015}, \cite{konstantinides:li:2016},
\cite{konstantinides:liu:passalidis:2025}, \cite{konstantinides:passalidis:2024j}, \cite{li:2016}, and others. These models often assume a renewal counting process or quasi-renewal counting processes where the inter-arrival times are identically distributed, and weak dependence is incorporated. However, few studies have dealt with more complex risk models where the counting process does not follow these standard structures, allowing for more sophisticated dependence relationships between the claim vectors.

For multivariate risk models driven by non-renewal, or non-quasi-renewal processes, the literature offers only a few contributions. For instance, in many of these papers, additional assumptions have been considered, such as the finite moment conditions for the counting process associated with the Matuszewska indices of the claim vectors. As a consequence, many subexponential claim amount distributions are not considered. Other supplementary conditions consist of specific tail conditions for the counting process and the claim amount distribution; see \cite{chen:li:cheng:2023}, \cite{yang:su:2023} and \cite{passalidis:2025} for more details.

Motivated by these limitations, in this paper we propose two general conditions (Assumption \ref{ass.KPX.3.2} for the finite-time horizon and Assumption \ref{ass.KPX.3.1} for the infinite-time horizon) that generalize the renewal framework in several directions. In particular, our conditions encompass quasi-renewal processes as well as inhomogeneous renewal processes (that is, processes with independent but not identically distributed inter-arrival times). In addition, our conditions cover many situations in which the inter-arrival times in the counting process are neither identically distributed nor independent. Such a general framework is particularly relevant in actuarial practice. For example, inhomogeneous renewal processes (such as inhomogeneous Poisson processes) can be used to model processes with seasonal patterns. For instance, in fire insurance, the number of claims tends to increase during the summer months.

Another situation in which the renewal process seems to defect in the risk model of \eqref{eq.KPX.1.2} can be described as follows. Because the $d$-lines of business share a common counting process, the occurrence of a claim in one particular insurance portfolio at time point $i$, say $X_1^{(i)}$, may increase the probability of future claim arrival in the other portfolios, namely the appearance of $X_k^{(j)}$, for $j > i$ and $k=2,\,\ldots,\,d$.  In order to better understand this scenario, let us assume that an insurance company has two insurance portfolios: one for car insurance and another for health insurance. Suppose that we are involved in a serious car accident and we have to pay for damage caused to the car, which is represented by $X_1^{(1)}$. After that, we could be asked to pay for some injuries sustained during the accident in the portfolio of health insurance and generate a sequence of claims such as $X_2^{(2)}\,,\;X_2^{(3)}\,,\;\ldots$. This scenario shows that claim arrival times could become smaller, and consequently there is a 'negative dependence' among the inter-arrival times. Such examples highlight the need for a more general framework allowing non-independent and non-identically distributed inter-arrival times.

We now proceed to address point (2). By interdependent claim vectors, we mean that there can be both component-wise as well as vector-wise dependencies.
In Section 2, the former can be characterized via multivariate heavy-tailed distribution classes, which include the range from independence to asymptotic dependence among the components of the claim vectors. For the latter, we use weak and sufficiently general dependencies, which include independence as a special case.

From one hand side, there is a need for component-wise dependencies so as to capture the interactions between the different portfolios of the insurer, while, on the other hand, dependencies among claim vectors play an important role for the evaluation of the total amount of the losses occurring at each claim arrival time $\tau_i.$ The example of the car insurance and health insurance policies discussed above is again relevant, but now attention is focused on the sizes of the claims rather than their arrival times.

As far as we know, only a few papers that consider continuous time risk models  with interdependent claims, namely \cite{chen:cheng:zheng:2025}, \cite{chen:konstantinides:passalidis:2025}, \cite{chen:wang:wang:2013}, \cite{gao:yang:2014},  \cite{konstantinides:passalidis:2025a}, and \cite{konstantinides:passalidis:2024j}. However, these works focus exclusively on renewal risk models.

The rest of the paper is organized as follows. In Section 2, we provide the preliminary concepts regarding the multivariate heavy-tailed distribution classes and the dependence structures, which will be employed. In Section 3, we give the two main results, for finite and the infinite-time horizons, together with their proofs and some auxiliary lemmas. Finally in Section 4, we present an application to ruin probabilities for finite and infinite-time horizons in a risk model with a general premium process and Brownian perturbations.

\section{Preliminaries} \label{sec.KPX.2}

Hereafter, all limit relations are understood as $\xto$. All vectors are denoted in bold, and their dimension is $d \in \bbn$. For two vectors ${\bf x}$, ${\bf z}$,
and some positive scalar $y>0$, we use the following standard notation
\beao
{\bf x}\pm{\bf z}= (x_1\pm z_1,\,\ldots,\,x_d\pm z_d)^{\top}\,,\; {\bf x}\odot{\bf z}= (x_1\,z_1,\,\ldots,\,x_d\,z_d)^{\top}\,, \; y\,{\bf x}=(y\,x_1,\,\ldots,\,y\,x_d)^{\top}\,.
\eeao
The vector ${\bf 0}=(0,\,\ldots,\,0)^{\top}$ denotes the origin of the axes. For two real numbers $a,\,b$, we write $a\vee b := \max\{a,\,b\}$ and $a\wedge b := \min\{a,\,b\}$. For any set $\bbb \in \overline{\bbr}^d := [-\infty,\,\infty]^d$, we denote by $\bbb^c$ its complement set, by $\overline{\bbb}$ its closed hull, by $\partial \bbb$ its border, and by ${\bf 1}_{\bbb}$ its indicator function.

Recall that a set $\bbb$ is called increasing if for any ${\bf x} \in \bbb$ and ${\bf z} \in \bbr_+^d$, it holds that ${\bf x}+{\bf z}\in \bbb$. For two positive uni-variate functions $f,\,g$, we write $f(x) \sim c\,g(x)$, for some $c \in (0,\,\infty)$, if
\beao
\lim \dfrac{f(x)}{g(x)} =c\,,
\eeao
we write $f(x) = o[g(x)]$, if
\beao
\lim \dfrac{f(x)}{g(x)} =0\,,
\eeao
write $f(x)=O[g(x)]$, if
\beao
\limsup \dfrac{f(x)}{g(x)} < \infty\,,
\eeao
and  $f(x) \asymp g(x)$ holds, if both $f(x)=O[g(x)]$ and $g(x)=O[f(x)]$, are true.

Correspondingly, if $\widetilde{f},\,\widetilde{g}$ are $d$-variate functions, and $\bbb \in \bbr^d \setminus \{{\bf 0}\}:= (-\infty,\,\infty)^d \setminus \{{\bf 0}\}$,
 we use the same asymptotic notations, for example we write $\widetilde{f}(x\,\bbb) =O[\widetilde{g}(x\,\bbb)]$, if
\beao
\limsup \dfrac{\widetilde{f}(x\,\bbb)}{\widetilde{g}(x\,\bbb)} < \infty\,.
\eeao
Finally for a random variable (or vector) $\Theta$, we write $\Theta \stackrel{d}{\sim} V$ to indicate that $\Theta$ has distribution $V$. The tail distribution of $V$ is denoted by $\bV(x) =1- V(x)$, for any $x \in \bbr$.

\subsection{Heavy tailed random variables} \label{sec.KPX.2.1}

Let us recall a one-dimensional distribution $V$ is said to be heavy-tailed, if $\bV(x) >0$ for all $x \in \bbr$, and
\beao
\int_{-\infty}^{\infty} e^{\vep \,x}\,V(dx) = \infty\,,
\eeao
for any $\vep >0$. Next, we introduce the following family of sets,
\beam \label{eq.KPX.2.1}
\mathscr{R} = \left\{A \subset \bbr^d\;:\;A\;\text{ open, increasing}\,,\;A^c\;\text{ convex}\,,\;{\bf 0}\notin \overline{A} \right\}\,,
\eeam
in order to define the rest subclasses of multivariate, heavy tailed distributions. Two important sets, from the family $\mathscr{R} $, which play a crucial role in actuarial practice, are the following:
\beam \label{eq.KPX.2.2}
A_{1} = \left\{{\bf z}\;:\;\sum_{i=1}^d l_i\,z_i > c \right\}\,,
\eeam
with $l_1,\,\ldots,\,l_d \geq 0$, $l_1+ \cdots + l_d =1$, and some $c>0$, and
\beam \label{eq.KPX.2.3}
A_{2} = \left\{{\bf z}\;:\;z_i > c_i\,,\;\exists \;i=1,\,\ldots,\,d \right\}\,,
\eeam
with $c_1,\,\ldots,\,c_d > 0$. In the special one-dimensional case, the set $A_2$ reduces to the interval $(c,\,\infty)$. See Remark \ref{rem.KPX.3.1} for the role of the sets $A_1$ and $A_2$ in actuarial practice.

Let ${\bf X} \stackrel{d}{\sim} F$ be a random vector, with distribution $F$, defined on $\bbr_+^d$. As shown in \cite{samorodnitsky:sun:2016}, for any $A \in \mathscr{R}$, the random variable
\beam \label{eq.KPX.2.4}
X_{A} := \sup \left\{ u\;:\;{\bf X} \in u\,A \right\}\,,
\eeam
follows a proper distribution $F_A$, whose tail satisfies
\beam \label{eq.KPX.2.5}
\bF_{A}(x) = \PP \left( \sup_{{\bf p} \in I_A} {\bf p}^{\top}\,{\bf X} >x \right) = \PP({\bf X} \in x\,A)\,,
\eeam
for any $x>0$, where $I_A$ is a suitable set of vectors whose existence for any $A \in \mathscr{R}$ was shown in \cite[Lemma 4.3(c)]{samorodnitsky:sun:2016}. Thus, applying \eqref{eq.KPX.2.1},  \eqref{eq.KPX.2.4} and  \eqref{eq.KPX.2.5}, the class of multivariate subexponential distributions was defined in \cite{samorodnitsky:sun:2016} as follows.

Let  $A \in \mathscr{R}$. We say that $F$ belongs to the class of multivariate subexponential distributions on $A$, symbolically $F \in \mathcal{S}_A$, if $F_A \in \mathcal{S} $, namely for any (or, equivalently, for some) integer $n \geq 2$ it holds
\beao
\lim \dfrac{\overline{F_A^{n*}}(x)}{\bF_A(x)}= n\,,
\eeao
where $F_A^{n*}$ represents the $n$-th order convolution power of $F_A$. Further, we
write
\beao
\mathcal{S}_\mathscr{R} := \bigcap_{A \in \mathscr{R}} \mathcal{S}_A\,,
\eeao
for the class of multivariate subexponential distributions.

Following the same line, the following two classes were introduced in \cite{konstantinides:passalidis:2024g}. Let  $A \in \mathscr{R}$. We say that $F$ belongs to the class of multivariate long tailed distributions on $A$, symbolically $F \in \mathcal{L}_A$, if $F_A \in \mathcal{L}$, namely for any (or, equivalently, for some) $y >0$ it holds
\beao
\lim \dfrac{\bF_A(x-y)}{\bF_A(x)}= 1\,.
\eeao
We say that $F$ belongs to the class of multivariate consistently varying distributions on $A$, symbolically $F \in \mathcal{C}_A$, if $F_A \in \mathcal{C}$,
namely it holds
\beao
\lim_{b\uparrow 1} \limsup \dfrac{\bF_A(b\,x)}{\bF_A(x)}= 1\,.
\eeao

Further, another general
multivariate distribution class was introduced in \cite{konstantinides:passalidis:2025h}, which also contains light-tailed distributions.

Let  $A \in \mathscr{R}$. We say that $F$ belongs to class of multivariate positively decreasing distributions on $A$, symbolically $F \in (\mathcal{P_D})_A$, if $F_A \in \mathcal{P_D}$, namely for any (or, equivalently, for some) $v >1$ it holds
\beao
\limsup \dfrac{\bF_A(v\,x)}{\bF_A(x)} < 1\,.
\eeao
We also denote $F \in (\mathcal{C} \cap \mathcal{P_D})_A$, if it holds $F_A \in \mathcal{C}\cap \mathcal{P_D}$. For all previous classes we keep the notation
\beao
\mathcal{B}_\mathscr{R} := \bigcap_{A \in \mathscr{R}} \mathcal{B}_A\,,
\eeao
for any $\mathcal{B} \in \{\mathcal{C} \cap \mathcal{P_D}\,,\;\mathcal{C}\,,\; \mathcal{P_D}\,,\;\mathcal{L}\}$.

Here, we deal mostly with the classes $\mathcal{S}_A$ and $(\mathcal{C} \cap \mathcal{P_D})_A$, for the description of the distribution of claim vectors. We should mention that the class $\mathcal{P_D}$ is general enough, and hence the class  $(\mathcal{C} \cap \mathcal{P_D})_A$ is only negligibly smaller than $\mathcal{C}_A$.

These classes have found applications in several multivariate risk models, see in
\cite{samorodnitsky:sun:2016}, \cite{chen:konstantinides:passalidis:2025},
\cite{konstantinides:liu:passalidis:2025}. For examples of distributions, that belong to classes $(\mathcal{C} \cap \mathcal{P_D})_A$ and $\mathcal{S}_A$, we refer to \cite[Sec. 4]{samorodnitsky:sun:2016} and \cite[Sec. 4]{konstantinides:liu:passalidis:2025}.

A famous class of heavy-tailed distributions in multidimensional set up is the (standard) multivariate regularly varying class, symbolically $MRV$. We remind that an one-dimensional distribution $V$ belongs to the class of regularly varying distributions, with index $\alpha \in (0,\,\infty)$, symbolically $V \in \mathcal{R}_{-\alpha}$, if for any $t>0$ it holds
\beao
\lim \dfrac{\bV(t\,x)}{\bV(x)} = t^{-\alpha}\,.
\eeao
We say that distribution $F$ belongs to $MRV$, if there exists some $V  \in \mathcal{R}_{-\alpha}$, with  $\alpha \in (0,\,\infty)$, and some non-degenerate to zero Radon measure $\mu$, such that
\beam \label{eq.KPX.2.6}
\lim \dfrac 1{\bV(x)} \PP \left( {\bf X} \in x\,\bbb \right) = \mu(\bbb)\,,
\eeam
for any Borel set $\bbb \subseteq \overline{\bbr}^d \setminus \{ {\bf 0}\}$, which is such that $\mu(\partial \bbb)=0$. Then we denote $F \in MRV(\alpha,\,\mu)$. The Radon measure $\mu$, possess the property of positive homogeneity, namely for any real $y>0$,
and Borel set $\bbb \subseteq \overline{\bbr}^d \setminus \{{\bf 0}\}$, it holds
\beam \label{eq.KPX.2.7}
\mu(y\,\bbb) = y^{-\alpha} \,\mu(\bbb)\,.
\eeam
For further reading about $MRV$ distributions, we refer the reader to \cite{resnick:2007}.

If $MRV$ represents the union of all $MRV(\alpha,\,\mu)$, for any  $\alpha \in (0,\,\infty)$ and any Radon measure $\mu$, then by combination of \cite[Proposition 2.1]{konstantinides:passalidis:2024g} and \cite[Proposition 3.1]{konstantinides:passalidis:2025h} we find the following inclusions:
\beam \label{eq.KPX.2.8}
MRV \subsetneq (\mathcal{C} \cap \mathcal{P_D})_\mathscr{R} \subsetneq \mathcal{C}_\mathscr{R} \subsetneq \mathcal{S}_\mathscr{R} \subsetneq \mathcal{L}_\mathscr{R}\,.
\eeam
These inclusions in \eqref{eq.KPX.2.8} remain valid if the classes $\mathcal{B}_\mathscr{R}$ are replaced by $\mathcal{B}_A$, for any $A \in \mathscr{R}$.
Consequently, the ordering of these classes in the one-dimensional case also extends to the multidimensional case. For further discussion on the ordering of heavy-tailed distributions in one dimension, we refer to \cite[Chapter 2]{leipus:siaulys:konstantinides:2023}.

Finally we recall the Matuszewska indexes, which are useful for characterizing the heavy tailed and related distributions classes. For an one-dimensional distribution $G$, with $\bG(x) > 0$, for all $x \in \bbr$, the upper and lower Matuszewska indexes are correspondingly defined as
\beao
J_G^+:= -\lim_{v\to \infty} \dfrac{\log \bG_*(v)}{\log v}\,,\qquad J_G^-:= -\lim_{v\to \infty} \dfrac{\log \bG^*(v)}{\log v}\,,
\eeao
where
\beao
\bG_*(v)=\liminf \dfrac{\bG(v\,x)}{\bG(x)}\,,\qquad \bG^*(v)=\limsup \dfrac{\bG(v\,x)}{\bG(x)}\,.
\eeao
Note that $G \in \mathcal{P_D}$ if and only if $J_G^- >0$.
If $G \in \mathcal{C}$, then $J_G^+ <\infty$, and if $G \in \mathcal{R}_{-\alpha}$,
then $J_G^-=J_G^+=\alpha$. For further discussions about the Matuszewska indexes we
refer to \cite[Subsection 2.1.2]{bingham:goldie:teugels:1987} and to \cite[subsection 2.4]{leipus:siaulys:konstantinides:2023}.

\subsection{Discussion on the dependence} \label{sec.KPX.2.2}

Now we introduce two main dependence structures that will be used to model the dependence among claim vectors. We first define their one-dimensional forms and then their multidimensional forms by reducing the problem to the one-dimensional form via the set $A \in \mathscr{R}$. The first dependence structure was introduced in \cite{lehmann:1966}, and since then it has been used in several topics related with subexponential distributions, see e.g. \cite{geluk:tang:2009}, \cite{geng:liu:wang:2023}, \cite{jiang:gao:wang:2014} among others. In order to avoid unnecessary details, we only define these structures for the nonnegative case.

\bde \label{def.KPX.2.1}
We say that the sequence of non-negative, random variables $\{Z_{i},\,i \in \bbn\}$ with corresponding distributions $\{G_{i},\,i \in \bbn\}$ are regression dependent, symbolically $RD$, if there exists some $x_0 > 0$, and some $K>0$, such that the inequality
\beao
\PP(Z_i> x_i\;\big|\;Z_j=x_j\,,\;j \in J_i) \leq K\,\bG_i(x_i)\,,
\eeao
holds for any $i \in \bbn$, $j \in J_i \subset \bbn \setminus \{i\}$, with $J_i \neq \emptyset$ and $x_i \wedge x_j > x_0$.
\ede

The $RD$ structure is a dependence structure, which is commonly employed when $G_i \in \mathcal{S}$, 
for the extension of results that were initially proved for independent random variables, see Remark \ref{rem.KPX.2.2} below. It is obvious that the $RD$ structure contains
the independence structure as a special case. However, in the case, where $G_i \in \mathcal{C} \subsetneq \mathcal{S}$, we can generalize statements on independent random variables via a more general dependence structure, namely the quasi asymptotic independence, denoted by $QAI$, which was introduced in \cite{chen:yuen:2009}. This dependence structure has been extensively employed in the literature, see e.g. \cite{chen:liu:2022}, \cite{cheng:2014}, \cite{dirma:nakiliuda:siaulys:2023}, \cite{li:2013}, among others.

\bde \label{def.KPX.2.2}
We say that the sequence of non-negative, random variables $\{Z_{i},\,i \in \bbn\}$ with distributions $\{G_{i},\,i \in \bbn\}$ are $QAI$, if it holds
\beao
\lim \dfrac{\PP(Z_i> x\,,\;Z_j> x)}{\bG_i(x)+ \bG_j(x)} =0\,,
\eeao
for any $i,\,j \in \bbn$, with  $i \neq j$.
\ede

The following two definitions represent the multivariate versions of Definitions \ref{def.KPX.2.1} and \ref{def.KPX.2.2}, appearing in \cite{konstantinides:passalidis:2024g}, and further used in \cite{chen:konstantinides:passalidis:2025}, for the description of interdependent claim vectors in a renewal driven risk model. For the sequence of non-negative, random vectors $\{{\bf Z}^{(i)},\,i \in \bbn\}$, we define
\beao
Z_A^{(i)} = \sup\left\{u\;:\;{\bf Z}^{(i)} \in u\,A \right\}\,,
\eeao
for any $i \in \bbn$.

\bde \label{def.KPX.2.3}
Let  $A \in \mathscr{R}$ be some fixed set. We say that the sequence of random vectors $\left\{{\bf Z}^{(i)},\,i \in \bbn \right\}$ is regression dependent on $A$, symbolically $RD_A$, if the sequence $\left\{ Z_A^{(i)},\,i \in \bbn \right\}$ is $RD$.
\ede

\bde \label{def.KPX.2.4}
Let  $A \in \mathscr{R}$ be some fixed set. We say that the sequence of random vectors $\left\{{\bf Z}^{(i)},\,i \in \bbn \right\}$ is quasi asymptotically independent on $A$, symbolically $QAI_A$, if the sequence $\left\{ Z_A^{(i)},\,i \in \bbn \right\}$ is $QAI$.
\ede

It is easy to see that for any $A \in \mathscr{R}$, it holds $RD_A \subsetneq QAI_A$. For examples of dependencies in Definitions \ref{def.KPX.2.3} and \ref{def.KPX.2.4} see in \cite[Subsection 2.2]{chen:konstantinides:passalidis:2025}.

The following remarks are mainly meant to provide information regarding the dependence structures that preserve the 
important properties of subexponential distributions, such as the single big jump principle, in one-dimensional as well as multidimensional cases.

\bre \label{rem.KPX.2.2}
A classical question of significant interest in applied probability concerns the single big jump principle. Specifically, suppose that we have random variables $Z_A^{(1)},\,\ldots,\,Z_A^{(n)}$ with distributions $G_A^{(1)},\,\ldots,\,G_A^{(n)}$, respectively, under what assumptions for the distribution class of $\left\{G_A^{(i)}\,,\;i=1,\,\ldots,\,n \right\}$ and the dependence structures for the  $\left\{Z_A^{(i)}\,,\;i=1,\,\ldots,\,n \right\}$ we obtain the asymptotic relation
\beam \label{eq.KPX.2.11}
\PP\left(\sum_{i=1}^n Z_A^{(i)}> x\right) \sim \sum_{i=1}^n\PP\left( Z_A^{(i)}> x\right) \,.
\eeam
This problem either in the form \eqref{eq.KPX.2.11}, or in more general weighted form, has been investigated in many papers, see in \cite{cheng:2014}, \cite{li:2013}, \cite{wang:2011}, \cite{yang:wang:leipus:siaulys:2013}, among many others.

In general, the following intuitive interpretation is true:
\begin{enumerate}
\item[(A)]
The smaller is the the distribution class of $\left\{G_A^{(i)}\,,\;i=1,\,\ldots,\,n \right\}$, the larger (weak) dependence structures can be used, in order to keep \eqref{eq.KPX.2.11} true.
\end{enumerate}

Although in the cases where $G_A^{(i)} \in \mathcal{S}$, with $J_{G_A^{(i)}}^+ <\infty$ for any $i=1,\,\ldots,\,n$ the dependence structures are 'well' ordered with respect to distribution classes (in order to keep \eqref{eq.KPX.2.11} true),
in the cases where $G_A^{(i)} \in \mathcal{S}$, with $J_{G_A^{(i)}}^+ =\infty$ for any $i=1,\,\ldots,\,n$ this is no more true. The reason for this, has to do with the two different approaches on this issue.

The first approach, which we shall follow, is based on the structure $RD$ from Definition \ref{def.KPX.2.1} and the structure introduced in \cite{ko:tang:2008}. These dependence structures allow comparison with the others in question (A) (as for example $QAI$ structure with $\mathcal{C}$ class), while in the same time these two assumptions can operate in complementary mode, in relation to question \eqref{eq.KPX.2.11}, see in \cite{jiang:gao:wang:2014} for more discussions on these two structures.

The second approach, given in \cite{foss:richards:2010}, considered that the $\{Z_A^{(i)}\,,\;i=1,\,\ldots,\,n\}$ are conditionally independent, given some $\sigma$-algebra. This idea undoubtedly has great theoretical value and appears general enough. However, in that paper, additional conditions are needed to establish \eqref{eq.KPX.2.11}, which seems rather restrictive. Most of their examples satisfying these conditions can be reduced to subcases of the first approach. Namely, in \cite[Example 1]{foss:richards:2010}, we can easily see that $G_i \in \mathcal{R}_{-\alpha}$, for some $\alpha \in (0,\,\infty)$, while the structure of dependence for the $Z_A^{(i)}$ is the $QAI$. In \cite[Example 4]{foss:richards:2010}, the $Z_A^{(i)}=\xi_i\,\eta_1,\,\eta_2,\,\cdots,\,\eta_i$, with $\{\xi_i\}$ independent and identically distributed (i.i.d.) random variables with distribution $F_{\xi} \in \mathcal{C}$, and with $\{\eta_i\}$ i.i.d. with distribution $F_{\eta}$, which belongs to the class of rapidly varying distributions, where $\{\xi_i \}$ independent of $\{\eta_i\}$. We can easily find out, through \cite[Lemma 3.2]{chen:yuan:2017}, that $G_A^{(i)} \in \mathcal{C}$ while the structure of dependence for the $Z_A^{(i)}A^{(i)}$ is the $QAI$. Even more, we can further generalize, with $\{ (\xi_i\,,\;\eta_i)\,,\; i \in \bbn\}$ sequence of i.i.d. pairs, $\xi_i,\,\eta_i$ arbitrarily dependent $F_{\xi} \in \mathcal{S}$, and with $J_{F_{\xi}}^+ < \infty$, and $\E[\eta^p]< \infty$, for $p > J_{F_{\xi}}^+$, and again from \cite[Lemma 3.2]{chen:yuan:2017} and \cite[Theorem 3.1]{geluk:tang:2009}, the relation \eqref{eq.KPX.2.11} remains true. Finally, \cite[Example 5]{foss:richards:2010} focuses in case, where $G_{A}^{(i)} \in \mathcal{S}$ with $J_{G_{A}^{(i)}}^+ = \infty$, but the dependence structure, that follows from multivariate normal distribution, represents a special case of $RD$, as we find in \cite[p. 876]{geluk:tang:2009}.

Therefore, we consider that in case, where $G_i \in \mathcal{S}$, it is plausible to use $RD$, rather than the conditions in \cite{foss:richards:2010} (see also in \cite{cheng:xu:2020}), although they are not overlapping each other.
\ere

\bre \label{rem.KPX.2.1}
The corresponding multivariate analogue of relation \eqref{eq.KPX.2.11}, can take two forms: 1) the multivariate linear single big jump principle and  2) the multivariate non-linear single big jump principle, see further for 2) in \cite{konstantinides:passalidis:2025a}. Here we put emphasis in linear form, which is satisfied by the classes in relation \eqref{eq.KPX.2.8}, except
class $\mathcal{L}_\mathscr{R}$, for independent random vectors. Namely, the (A) is reduced to
\beam \label{eq.KPX.2.12}
\PP\left(\sum_{i=1}^n {\bf Z}^{(i)} \in  x\,A \right) \sim \sum_{i=1}^n\PP\left( {\bf Z}^{(i)} \in  x\,A\right) \,,
\eeam
for any $A \in \mathscr{R}$. The name 'linear' stems from the fact that relation \eqref{eq.KPX.2.12} is insensitive with respect to the dimension of the random vectors. Hence, due to this insensitivity with respect to dimensionality, the question of (A) for \eqref{eq.KPX.2.12}, can be easily transferred to the question of (A) for \eqref{eq.KPX.2.11}, by employing the corresponding multivariate distribution classes and dependence structures on set $A \in \mathscr{R}$.
\ere

\section{Main results} \label{sec.KPX.3}

Now we can give the two main results of the paper.

\subsection{Finite-time horizon} \label{sec.KPX.3.1}

In this section we assume that the discounted aggregate claims have claim vectors that permit the $RD_A$ structure, with distribution from class $\mathcal{S}_A$.

We also consider the following general assumption for the counting process $\{N(t)\,,\;t\geq 0\}$, that is related with its second factorial moment measure, symbolically $\alpha^{(2)}$, which is defined by
\beao
\alpha^{(2)}(ds,dt):=E\left[N(ds)\,N(dt)\right]-E[N(ds)\delta_s(dt)],
\eeao
where $\delta_s(dt)$ is the Dirac mass on the diagonal $t=s$. This means that the second factorial moment measure is the mean measure of $N(ds)N(dt)$ with the diagonal $t=s$ removed.

\begin{assumption} \label{ass.KPX.3.2}
Let consider a fixed $T \in \Lambda \setminus \{\infty\}$. We assume that there exists $C<\infty$ such that
\beao
\alpha^{(2)}(ds,dt) \le C\, m(ds)m(dt)\,,
\eeao
on $[0,\,T]^2$.
\end{assumption}

\bre \label{rem.KPX.3.a}
Note that for any inhomogeneous renewal process, because counts on disjoint sets are independent, it holds
\begin{eqnarray*}
\alpha^{(2)}(ds,dt)&=&E\left[N(ds)\,N(dt)\right]-E[N(ds)\delta_s(dt)]\nonumber\\
&=&E\left[N(ds)\right]E\left[N(dt)\right]+E[N(ds)\delta_s(dt)]-E[N(ds)\delta_s(dt)]\nonumber\\
&=&E\left[N(ds)\right]E\left[N(dt)\right]=m(ds)m(dt)\,.
\end{eqnarray*}
So Assumption \ref{ass.KPX.3.2} holds with $C=1$ for all $T \in \Lambda \setminus \{\infty\}$.
\ere

\bre \label{rem.KPX.3.b}
It is a fact that Assumption \ref{ass.KPX.3.2}, does not include only inhomogeneous renewal processes but it can cover also cases where the $\{N(t)\,,\;t\geq 0\}$, does not have necessarily neither independent, nor identically distributed increments. For example, it contains all the processes, such that satisfy
\beam \label{eq.KPX.3.a}
\E[N(T)\,N(T)] < \infty\,,
\eeam
for fixed $T \in \Lambda \setminus \{\infty\}$. Indeed,  let $\underline{t}:= \inf\{t\;:\;t \in \Lambda\}$. Let $\PP(\tau_1=\underline{t})>0$. We examine two subcases. Firstly, if either $s \in [0,\,\underline{t})$, or  $t \in [0,\,\underline{t})$, then it holds
\beao
\alpha^{(2)}(ds,dt) \leq \E[N(ds)\,N(dt)] =0\,,
\eeao
and therefore Assumption  \ref{ass.KPX.3.2} holds trivially, as $0=0$.

For the second case, if $s,\,t \in [\underline{t},\,T]$, then we obtain
\beao
\dfrac{\alpha^{(2)}(ds,dt)}{m(ds)\,m(dt)} \leq \dfrac{\E[N(ds)\,N(dt)] }{m(ds)\,m(dt)} \leq \dfrac{\E[N(T)\,N(T)] }{[m(\underline{t})]^2}\,,
\eeao
and from relation \eqref{eq.KPX.3.a}, the last ratio is bounded from above. In the case where $\PP(\tau_1=\underline{t})=0$, we examine the two subcases, for $s \in [0,\,\underline{t}]$ or $t \in [0,\,\underline{t}]$, and for $s,\,t \in (\underline{t},\,\infty]$, and further the methodology is similar. We can see that as smaller the $T \in \Lambda$, so easier relation \eqref{eq.KPX.3.a} is satisfied. This shows the generality - flexibility provided by Assumption  \ref{ass.KPX.3.2}.
\ere

\bth \label{th.KPX.3.1}
Let $A \in \mathscr{R}$ be some fixed set, and the discounted aggregate claims of \eqref{eq.KPX.1.2}. We suppose that Assumptions \ref{ass.KPX.1.1} and  \ref{ass.KPX.3.2} with fixed  $T \in \Lambda \setminus \{\infty\}$ hold, and $\left\{{\bf X}^{(i)},\,i \in \bbn \right\}$ are $RD_A$, with $F \in \mathcal{S}_A$. Then it holds
\beam \label{eq.KPX.3.1}
\PP\left({\bf D}(T) \in  x\,A \right) \sim \int_{0}^T \PP\left( {\bf X}\,e^{-r\,s} \in x\,A\right)\,m(ds) \,.
\eeam
\ethe

\bre \label{rem.KPX.3.1}
If in conditions of Theorem \ref{th.KPX.3.1} we require $F \in \mathcal{S}_\mathscr{R}$, and the structures of $\left\{{\bf X}^{(i)},\,i \in \bbn \right\}$ are $RD_A$, for any $A \in \mathscr{R} $, then we can imply, from the proof of Theorem \ref{th.KPX.3.1}, that relation \eqref{eq.KPX.3.1} holds for any $A \in \mathscr{R}$. Although relation \eqref{eq.KPX.3.1} eventually seems complicated in calculation, it can take several simpler expressions, when we determine the set $A$, the dependence structures of the marginals, and sometimes also the distribution class of marginals, see in \cite[Section 4]{konstantinides:liu:passalidis:2025}, for related discussions.

Further, the sets $A_1$ and $A_2$ from \eqref{eq.KPX.2.2} and \eqref{eq.KPX.2.3},
respectively, represent two important sets for the insurance industry, since provide two different kind of control for the insurer's solvency. Concretely, the entrance probability of the discounted aggregate claims into set $A_1$, indicates the probability of the event, the sum of the discounted aggregate claims of $d$-lines to exceed some 'high' threshold, while the corresponding entrance into $A_2$, shows the probability in one of the $d$-lines the discounted aggregate claims to exceed some 'high' threshold. This 'high' threshold usually has relation with the insurer's initial capital.

Let us note that in one-dimensional sub-case, where $A=(1,\,\infty)$, this result continuous to be new, and in some sense is similar to \cite[Theorem 1]{cheng:xu:2020}, where was considered the ruin probability over finite-time horizon, but there was employed conditional independent claims instead of $RD$, and a inhomogeneous Poisson driven risk model.
\ere

The following corollary presents a more explicit asymptotic expression, in comparison with \eqref{eq.KPX.3.1}, when the distribution $F$ is restricted in $MRV$. We obserbe that for any $A \in \mathscr{R}$, we obtain $\mu(A) \in (0,\,\infty)$, see in \cite[proof of Proposition 4.14]{samorodnitsky:sun:2016}, while if $F \in MRV(\alpha,\,\mu)$, then there exists $V \in R_{-\alpha}$, recall \eqref{eq.KPX.2.6}.

\bco \label{cor.KPX.3.1}
Under the conditions of Theorem \ref{th.KPX.3.1}, with restriction that $F \in MRV(\alpha,\,\mu)$, with $\alpha \in (0,\,\infty)$, then we obtain
\beao
\PP\left({\bf D}(T) \in  x\,A \right) \sim \mu(A)\,\bV(x)\,\int_{0}^T \,e^{-\alpha\,r\,s}\,m(ds) \,.
\eeao
\eco

\subsection{Infinite-time horizon} \label{sec.KPX.3.2}

For the infinite-time horizon we need some alternative conditions, which guarantee the finiteness of the sums, like that in \eqref{eq.KPX.1.2} with $T = \infty$, are non-defective, namely they have not mass on $\infty$.

Some sufficient conditions, are provided by the requirement of positive and finite Matuszewska indexes on the claim vectors, in combination with the fact that $r>0$, and the following moment condition.

\begin{assumption} \label{ass.KPX.3.1}
Let $A \in \mathscr{R}$ be a fixed set, $r>0$,  and $F \in (\mathcal{C} \cap \mathcal{P_D})_A$. We assume that there exist some $q_1,\,q_2$, with $0<q_1 < J_{F_A}^- \leq J_{F_A}^+ <q_2 < \infty$, such that it holds
\beam \label{eq.KPX.3.3}
\sum_{i=1}^{\infty} \left(\E\left[e^{-q_1\,r\,\tau_i}\right] \bigvee \E\left[e^{-q_2\,r\,\tau_i}\right] \right)^{1/q} <\infty\,,
\eeam
where $q ={\bf 1}_{\left\{0 <J_{F_A}^+ <1 \right\}} +q_2\,{\bf 1}_{\left\{J_{F_A}^+ \geq 1 \right\}}$.
\end{assumption}

Although Assumption \ref{ass.KPX.3.1} is more restrictive than Assumption \ref{ass.KPX.3.2}, used for the finite time horizon, it continuous to remain quite general. The following examples demonstrate this generality.

\bexam \label{exam.KPX.3.1}
We consider that the inter-arrival times $\{\theta_i:=\tau_i -\tau_{i-1}\,,\;i \in \bbn\}$ are widely lower orthant dependent, symbolically $WLOD$, namely there exists a sequence of finite positive numbers $\{g_L(n)\,,\;n \in \bbn\}$, such that for any $x_1,\,\ldots,\,x_n$ to hold
\beam \label{eq.KPX.ex.3.a}
\PP\left( \bigcap_{i=1}^n \{\theta_i \leq x_i\}\right) \leq g_L(n) \prod_{i=1}^n \PP(\theta_i \leq x_i)\,,
\eeam
see \cite{wang:wang:gao:2013} for the introduction of $WLOD$.

Further, we suppose that the $\{\theta_i\,,\;i \in \bbn\}$ are identically distributed. Hence, we obtain that the $\{N(t)\,,\;t\geq 0\}$ is a quasi-renewal process. If additionally there exists some $\delta \in (0,\,-\log \E[e^{-r\,J_{F_A}^-\,\theta_1}])$, such that it holds
\beam \label{eq.KPX.ex.3.b}
\limsup_{\nto} g_L(n) \,e^{-\delta\,n} < \infty\,,
\eeam
then Assumption \ref{ass.KPX.3.1} holds. This can be verified following the lines of the proof of \cite[Th. 2]{yang:wang:2012}. We note that relation \eqref{eq.KPX.ex.3.a} represents a general enough dependence structure and the additional condition \eqref{eq.KPX.ex.3.b} is quite mild. We refer to \cite{wang:wang:gao:2013} and \cite{wang:cui:wang:ma:2012} for more details about \eqref{eq.KPX.ex.3.a} and \eqref{eq.KPX.ex.3.b}.
\eexam

It is obvious that if the counting process $\{N(t)\,,\; t\geq 0\}$ represents a renewal process, then relation \eqref{eq.KPX.ex.3.b} is true for all $\delta>0$ and hence also Assumption \ref{ass.KPX.3.1} via previous example.

\bexam \label{exam.KPX.3.2}
Let consider the inter-arrival times  $\{\theta_i\,,\;i \in \bbn\}$, that are arbitrarily dependent with distributions, which are bounded from below, namely it holds $\theta_i \geq a_i$, almost surely, with $a_i >0$. Such cases are very realistic, especially in life insurance. Let $\hat{a} = \inf_{n\geq 1} \left( \bigwedge_{i=1}^n a_i\right)$. Then it holds
\beao
&&\sum_{i=1}^{\infty} \left(\E\left[e^{-q_1\,r\,\tau_i} \right]\bigvee \E\left[e^{-q_2\,r\,\tau_i} \right] \right)^{1/q} \leq \sum_{i=1}^{\infty} \left(\E\left[e^{-q_1\,r\,\sum_{j=1}^i \theta_j} \right] + \E\left[e^{-q_2\,r\,\sum_{j=1}^i \theta_j} \right] \right)^{1/q}\\[2mm]
&&\leq \sum_{i=1}^{\infty} \left(\left[e^{-q_1\,r\,\hat{a}} \right]^i + \left[e^{-q_2\,r\,\check{a}} \right]^i \right)^{1/q} \leq \sum_{i=1}^{\infty} \left[e^{-q_1\,r\,\check{a}} \right]^{i/q} + \sum_{i=1}^{\infty} \left[e^{-q_2\,r\,\check{a}} \right]^{i/q} < \infty\,,
\eeao
where at the pre-last step we used the elementary inequality $(x+y)^z \leq x^z + y^z$, for any $x,\,y \geq 0$, and any $z \in (0,\,1]$.
\eexam

In spite of the restrictions of Assumption \ref{ass.KPX.3.1}, in comparison with the assumptions of Theorem \ref{th.KPX.3.1}, namely \eqref{eq.KPX.3.3} and $F \in (\mathcal{C} \cap \mathcal{P_D})_A $,  the following result permits dependence structure $QAI_A \supsetneq RD_A $, for the claim vectors.

\bth \label{th.KPX.3.2}
Let $A \in \mathscr{R}$ be some fixed set, and consider the discounted aggregate claims of relation \eqref{eq.KPX.1.2}. We suppose that Assumptions \ref{ass.KPX.1.1} and \ref{ass.KPX.3.1} are valid, and the dependence structure of $\left\{{\bf X}^{(i)},\,i \in \bbn \right\}$ is $QAI_A$. Then we find
\beam \label{eq.KPX.3.4}
\PP\left({\bf D}(\infty) \in  x\,A \right) \sim \int_{0}^{\infty} \PP\left( {\bf X}\,e^{-r\,s} \in  x\,A\right)\,m(ds)\,.
\eeam
\ethe

Similarly with what was noticed after  Theorem \ref{th.KPX.3.1}, here also, if we consider that $F \in (\mathcal{C} \cap \mathcal{P_D})_\mathscr{R}$, $\left\{{\bf X}^{(i)},\,i \in \bbn \right\}$ are $QAI_A$, for any $A \in \mathscr{R}$, and there exist $q_1(A),\,q_2(A)$, for any $A \in \mathscr{R}$, such that relation \eqref{eq.KPX.3.3} holds, then \eqref{eq.KPX.3.4} is true for any  $A \in \mathscr{R}$.

The following corollary gives a more explicit form of \eqref{eq.KPX.3.4} under the
restriction that $F \in MRV$. The proof of this corollary, follows the same argument as in proof of Corollary \ref{cor.KPX.3.1}, and therefore is omitted. We should mention that if $F \in MRV(\alpha,\,\mu)$, then $F_A \in \mathcal{R}_{-\alpha}$, for any  $A \in \mathscr{R}$, and hence $J_{F_A}^- = J_{F_A}^+ = \alpha$.

\bco \label{cor.KPX.3.2}
Under the conditions of  Theorem \ref{th.KPX.3.2}, under the additional restriction that
$F \in MRV(\alpha,\,\mu)$, with $\alpha \in (0,\,\infty)$, we have the relation
\beao
\PP\left({\bf D}(\infty) \in  x\,A \right) \sim \mu(A)\,\bV(x)\,\int_{0}^{\infty} \,e^{-\alpha\,r\,s}\,\,m(ds) \,.
\eeao
\eco
We note that the results of Theorem \ref{th.KPX.3.2} and Corollary \ref{cor.KPX.3.2} are new even for the one-dimensional subcase, with $A = (1,\,\infty)$.

\bre \label{rem.KPX.b}
We should mention that the proof of Theorem \ref{th.KPX.3.2} could not be derived by the current result from stochastic recurrence equations, as it happens several times for the ruin probability over infinite time horizon, see for example in \cite{konstantinides:mikosch:2005}. Indeed, our proof is based on the multivariate linear single big jump principle of some infinite randomly weighted sums. Although the framework of stochastic recurrence equations in multidimensional set up provides a more general model, due to random matrices as random weights, the current results can be applied to this model, for $T=\infty$, only in the case, where the $\left\{\left({\bf X}^{(i)},\,e^{-r\,\theta_i} \right)\,,\; i \in \bbn\right\}$ are i.i.d. copies of $\left({\bf X},\,e^{-r\,\theta_1} \right)$ and $F \in MRV(\alpha,\,\mu)$. Hence, it holds for renewal risk model with independent claim vectors, which we aim to generalize in the nonstandard model, see in $(1)$ and $(2)$ in Subsection \ref{subsection.1.2}. For more details about the stochastic recurrence equations see in \cite{buraczewski:damek:mikosch:2016}.
\ere

\subsection{Argumentation} \label{sec.KPX.3.3}

Before the proof of Theorem \ref{th.KPX.3.1}, we need two preliminary lemmas, that are related with the multivariate linear single big jump principle for (randomly) weighted sums.

\ble \label{lem.KPX.3.1}
Let $A \in \mathscr{R}$ be some fixed set, and choose some fixed $n \in \bbn$. We
consider that the ${\bf Z}^{(1)},\,\ldots,\,{\bf Z}^{(n)}$ are nonnegative $RD_A$ random vectors with distributions $G_1,\,\ldots,\,G_n \in \mathcal{S}_A$, respectively, and $G_i*G_j \in \mathcal{S}_A$, for any $1\leq i \neq j \leq n$. Then, for some constants $0<a \leq b < \infty$, it holds
\beam \label{eq.KPX.3.5-1}
\lim \sup_{{\bf c}_n \in [a,\,b]^n} \left| \dfrac{\PP\left(\sum_{i=1}^n c_i\,{\bf Z}^{(i)} \in  x\,A \right)}{\sum_{i=1}^n\,\PP\left( c_i\,{\bf Z}^{(i)} \in  x\,A \right)} - 1 \right|=0\,,
\eeam
where ${\bf c}_n:=(c_1,\,\ldots,\,c_n)$.
\ele

\pr~
At first, we notice that since $G_i,\,G_j \in \mathcal{S}_A$ and $G_i*G_j \in \mathcal{S}_A$, then by \cite[Th. 3.4]{konstantinides:passalidis:2024g} we obtain $G_A^{(i)}*G_A^{(j)} \in \mathcal{S}$. Hence, for the upper bound of \eqref{eq.KPX.3.5-1} it holds
\beam \label{eq.KPX.3.5-2} \notag
\PP\left(\sum_{i=1}^n c_i\,{\bf Z}^{(i)} \in  x\,A \right) &\leq& \PP\left(\sum_{i=1}^n c_i\,Z_A^{(i)} > x \right) \\[2mm]
&\sim& \sum_{i=1}^n\PP\left( c_i\,Z_A^{(i)} > x \right)=\sum_{i=1}^n\PP\left( c_i\,{\bf Z}^{(i)} \in  x\,A \right)\,,
\eeam
uniformly for ${\bf c}_n \in [a,\,b]^n$, where at the first step we used \cite[Prop. 2.4]{konstantinides:passalidis:2024g}, and at the second we applied \cite[Lem. 8]{wang:2011}.

For the lower bound of \eqref{eq.KPX.3.5-1}, since the summands are nonnegative and the $x\,A \in \mathscr{R}$ is an increasing set (recall that $\mathscr{R}$ represents a cone with respect to positive scalar multiplication), we find via Bonferroni's inequality that it holds
\beam \label{eq.KPX.3.5-3}
\PP\left(\sum_{i=1}^n c_i\,{\bf Z}^{(i)} \in  x\,A \right) &\geq& \PP\left(\bigcup_{i=1}^n \left\{ c_i\,{\bf Z}^{(i)} \in  x\,A \right\} \right) \\[2mm] \notag
&\geq& \sum_{i=1}^n \PP\left( c_i\,{\bf Z}^{(i)} \in  x\,A \right) - \sum_{1\leq i\neq j \leq n} \PP\left( c_i\,{\bf Z}^{(i)} \in  x\,A\,,\;c_j\,{\bf Z}^{(j)} \in  x\,A \right).
\eeam
for all ${\bf c}_n \in [a,\,b]^n$. For each term of the last sum in \eqref{eq.KPX.3.5-3}, since the ${\bf Z}^{(1)},\,\ldots,\,{\bf Z}^{(n)}$ are $RD_A$, and ${\bf c}_n \in [a,\,b]^n$, we can find some sufficiently large $x_0^*>0$, such that for any $x\geq x_0^*$ it holds
\beam \label{eq.KPX.3.5-4}
&&\PP\left( c_i\,{\bf Z}^{(i)} \in  x\,A\,,\;c_j\,{\bf Z}^{(j)} \in  x\,A \right)= \PP\left( c_i\,Z_A^{(i)} > x\,,\;c_j\, Z_A^{(j)} > x \right)\\[2mm] \notag
&&\leq K\,\PP\left( c_i\,Z_A^{(i)} > x \right)\,\PP\left( c_j\, Z_A^{(j)} > x \right) \leq K\,\PP\left( c_i\,{\bf Z}^{(i)} \in  x\,A \right)\,\PP\left( c_j\,{\bf Z}^{(j)} \in  x\,A \right).
\eeam
So, from relations \eqref{eq.KPX.3.5-3} and \eqref{eq.KPX.3.5-4} we obtain
\beam \label{eq.KPX.3.5-5}
\PP\left(\sum_{i=1}^n c_i\,{\bf Z}^{(i)} \in  x\,A \right) \gtrsim [1-o(1)]\,\sum_{i=1}^n\,\PP\left( c_i\,{\bf Z}^{(i)} \in  x\,A \right)\,,
\eeam
uniformly for ${\bf c}_n \in [a,\,b]^n$. From \eqref{eq.KPX.3.5-2}, \eqref{eq.KPX.3.5-5} in combination with their uniformity with respect to ${\bf c}_n \in [a,\,b]^n$, we conclude \eqref{eq.KPX.3.5-1}.
~\halmos

The next lemma clarifies the concept of multivariate linear single big jump principle for the weighted sum of Lemma \ref{lem.KPX.3.1}.

\ble \label{lem.KPX.3.2}
Let $A \in \mathscr{R}$ be some fixed set, and choose some fixed $n \in \bbn$. Under the assumptions of Lemma \ref{lem.KPX.3.1}, for some constants $0<a \leq b < \infty$, it holds
\beam \label{eq.KPX.3.5-6}
\lim \sup_{{\bf c}_n \in [a,\,b]^n} \dfrac{\PP\left(\sum_{i=1}^n c_i\,{\bf Z}^{(i)} \in  x\,A\,,\;\bigcap_{k=1}^n \left\{c_k\,{\bf Z}^{(k)} \notin x\,A \right\} \right)}{\sum_{i=1}^n\,\PP\left( c_i\,{\bf Z}^{(i)} \in  x\,A \right)} =0\,.
\eeam
\ele

\pr~
We have
\beam \label{eq.KPX.3.5-7}
&&\PP\left(\sum_{i=1}^n c_i\,{\bf Z}^{(i)} \in  x\,A\,,\;\bigcap_{k=1}^n \left\{c_k\,{\bf Z}^{(k)} \notin x\,A \right\} \right) \\[2mm] \notag
&&=\PP\left(\sum_{i=1}^n c_i\,{\bf Z}^{(i)} \in  x\,A \right) - \PP\left(\sum_{i=1}^n c_i\,{\bf Z}^{(i)} \in  x\,A\,,\;\bigcup_{k=1}^n \left\{c_k\,{\bf Z}^{(k)} \in x\,A \right\} \right) \\[2mm] \notag
&&=\PP\left(\sum_{i=1}^n c_i\,{\bf Z}^{(i)} \in  x\,A \right) - \PP\left(\bigcup_{k=1}^n \left\{c_k\,{\bf Z}^{(k)} \in x\,A \right\} \right) \leq \PP\left(\sum_{i=1}^n c_i\,{\bf Z}^{(i)} \in  x\,A \right)\\[2mm] \notag
&& - \sum_{k=1}^n\PP\left( c_k\,{\bf Z}^{(k)} \in x\,A \right) + \sum_{1\leq k <j \leq n} \PP\left( c_k\,{\bf Z}^{(k)} \in x\,A\,,\;c_j\,{\bf Z}^{(j)} \in x\,A \right)\\[2mm] \notag
&& \sim \sum_{1\leq k <j \leq n} \PP\left( c_k\,{\bf Z}^{(k)} \in x\,A\,,\;c_j\,{\bf Z}^{(j)} \in x\,A \right) = o\left[\sum_{i=1}^n \PP\left( c_i\,{\bf Z}^{(i)} \in  x\,A \right)\right]\,,
\eeam
uniformly for ${\bf c}_n \in [a,\,b]^n$, where at the second step we used the fact that the $x\,A$ is increasing and the $c_i\,{\bf Z}^{(i)}$, for $i=1,\,\ldots,\,n$, are nonnegative. At the third step, we used  Bonferroni's inequality, and at the fourth Lemma \ref{lem.KPX.3.1}. The last step follows through the same reasoning with relation \eqref{eq.KPX.3.5-4}.
~\halmos

\bre \label{rem.KPX.a.1}
The conditions of Lemmas \ref{lem.KPX.3.1}, and \ref{lem.KPX.3.2}, are more general than that used in the proof of Theorem \ref{th.KPX.3.1}, since it is enough the identical distribution of $G_i=G$, for $i=1,\,\ldots,\,n$. We observe also that in such a case $G^{2*} \in \mathcal{S}_A$, since
\beao
G^{2*}(x\,A) \sim 2\,G(x\,A)\,,
\eeao
see in \cite[Cor. 4.10, Prop. 4.12(a)]{samorodnitsky:sun:2016}. For more details about the closure properties of $\mathcal{S}_A$ with respect to convolution, we refer to \cite[Sec. 3.2]{konstantinides:passalidis:2024g}.
\ere

\bre \label{rem.KPX.a.2}
Under the assumptions of Lemmas \ref{lem.KPX.3.1}, and \ref{lem.KPX.3.2}, for the ${\bf Z}^{(1)},\,\ldots,\,{\bf Z}^{(n)}$ and if the $\Theta_1,\,\ldots,\,\Theta_n$ are arbitrarily dependent random variables, independent of the ${\bf Z}^{(1)},\,\ldots,\,{\bf Z}^{(n)}$, with each of them defined on $[a,\,b]$, with $0<a \leq b < \infty$, then by integration with respect to the $\Theta_1,\,\ldots,\,\Theta_n$, and via dominated convergence theorem we find that it holds
\beam \label{eq.KPX.3.5-8}
\PP\left(\sum_{i=1}^n \Theta_i\,{\bf Z}^{(i)} \in  x\,A \right) \sim \sum_{i=1}^n \PP\left( \Theta_i\,{\bf Z}^{(i)} \in  x\,A \right)\,,
\eeam
and
\beam \label{eq.KPX.3.5-9}
\PP\left(\sum_{i=1}^n \Theta_i\,{\bf Z}^{(i)} \in  x\,A\,,\;\bigcap_{k=1}^n \left\{\Theta_k\,{\bf Z}^{(k)} \notin x\,A \right\} \right)=o\left[ \sum_{i=1}^n \PP\left( \Theta_i\,{\bf Z}^{(i)} \in  x\,A \right)\right]\,.
\eeam
Further, if the conditions in Lemmas \ref{lem.KPX.3.1}, and \ref{lem.KPX.3.2}, are reduced to $G_1,\,\ldots,\,G_n \in \mathcal{S}_\mathscr{R}$, $G_i*G_j \in \mathcal{S}_\mathscr{R}$, for any $1\leq i \neq j \leq n$, and the ${\bf Z}^{(1)},\,\ldots,\,{\bf Z}^{(n)}$ are $RD_A$, for any $A \in \mathscr{R}$, then relations \eqref{eq.KPX.3.5-1}, \eqref{eq.KPX.3.5-6} (and similarly, for \eqref{eq.KPX.3.5-8}, \eqref{eq.KPX.3.5-9}) hold for any $A \in \mathscr{R}$.
\ere

\noindent{\bf Proof of Theorem \ref{th.KPX.3.1}}~
At first, for sake of simplicity of the notation we write
\beam \label{eq.KPX.3.5-10}
J_x(T):=\sum_{\tau_i\leq T} {\bf 1}_{\left\{ \,{\bf X}^{(i)}\,e^{-r\,\tau_i} \in\; x\,A \right\}}\,,\qquad \Lambda_x : =\E\left[ J_x(T)\right]\,.
\eeam
Further, we denote by $\mathcal{F}_{N}:=\sigma\left(N(\bbb)\;:\;\bbb \subseteq [0,\,T] \right)$, the $\sigma$-algebra, generated by the counting process $\{N(t)\,,\;t \in [0,\,T]\}$. The conditional expectation with respect to $\mathcal{F}_{N}$, means conditional with respect of all the moments $\{\tau_i\,,\; i \in \bbn\}$, since $N(T) = \sum_{i=1}^{\infty} {\bf 1}_{\{\tau_i \leq T\}}$. By $\E_N$ we denote the expectation with respect to $\{N(t)\,,\;t \in [0,\,T]\}$. Hence, for the $\Lambda_x$ from \eqref{eq.KPX.3.5-10}, by the conditional expectation with respect to  $\mathcal{F}_{N}$, we obtain
\beam \label{eq.KPX.3.5-11}
&&\Lambda_x=\E \left[ \E_N \left(J_x(T)\;\big|\;\mathcal{F}_{N} \right) \right] =\E \left[ \E_N \left(\sum_{\tau_i\in N([0,\,T])} {\bf 1}_{\left\{ {\bf X}^{(i)}\,e^{-r\,\tau_i}\, \in\; x\,A \right\}}\;\Big|\;\mathcal{F}_{N} \right) \right] \\[2mm] \notag
&&= \E \left[ \sum_{\tau_i\leq T}\E_N \left( {\bf 1}_{\left\{ {\bf X}^{(i)}\,e^{-r\,\tau_i}\, \in\; x\,A \right\}}\;\Big|\;\mathcal{F}_{N} \right) \right]= \E \left[ \sum_{\tau_i\leq T} \PP\left( {\bf X}^{(i)}\,e^{-r\,\tau_i} \in x\,A \;\Big|\;\mathcal{F}_{N} \right) \right] \\[2mm] \notag
&&=\E \left[ \sum_{\tau_i\leq T} \PP\left( {\bf X}^{(i)}\,e^{-r\,\tau_i} \in x\,A \right) \right] =\E \left[ \int_{0}^T \PP\left( {\bf X}\,e^{-r\,s} \in x\,A \right)\,N(ds) \right] \\[2mm] \notag
&&= \int_{0}^T \PP\left( {\bf X}\,e^{-r\,s} \in x\,A \right)\,m(ds)\,,
\eeam
where at the third step we used the linearity of expectation and Fubini's theorem, at the fifth step we took into account the property that $\PP\left(Y \in B(Z)\;\big|\;\mathcal{F} \right)= \PP\left(Y \in B(z)\right)\,\big|_{Z=z}$, for a $\mathcal{F}$-measurable random variable $Z$, and the $Y$ is independent of $\mathcal{F}$. The sixth step follows by definition of the integral of random measure, while at the last step we applied the Campell's theorem.

Hence, it is enough to show that
\beam \label{eq.KPX.3.5-12}
\PP\left( {\bf D}(T) \in x\,A \right)\sim \Lambda_x\,.
\eeam
Let us note, that it holds
\beam \label{eq.KPX.3.5-13}
\PP\left( J_x(T) =1 \right)\leq \PP \left( {\bf D}(T) \in x\,A \right) = \PP\left( J_x(T) \geq 1 \right) + \PP \left( {\bf D}(T) \in x\,A\,,\;J_x(T) =0 \right) \,.
\eeam
Having in mind, the single big jump principle from Lemmas \ref{lem.KPX.3.1}, and \ref{lem.KPX.3.2}, we can show that relation \eqref{eq.KPX.3.5-12} is true, through the following three steps.

\begin{enumerate}
\item
$\PP\left( J_x(T) \geq 2 \right) = o(\Lambda_x)$,
\item
$\PP\left( J_x(T) \geq 1 \right) \sim \Lambda_x$,
\item
$\PP \left( {\bf D}(T) \in x\,A\,,\;J_x(T) =0 \right)=o(\Lambda_x)$.
\end{enumerate}

We begin with the first step. We can see that it holds
\beao
{\bf 1}_{\left\{ J_x(T) \geq 2  \right\}} \leq \dfrac 12\,J_x(T)[J_x(T) -1]\,,
\eeao
almost surely (recall also that the set $x\,A$ is increasing). Taking expectations of both sides of the last inequality we obtain
\beam \label{eq.KPX.3.5-14}
\PP\left( J_x(T) \geq 2 \right)\leq  \dfrac 12\,\E\left[J_x(T)\,(J_x(T) -1) \right]\,.
\eeam
Due to the independence of the indicator function, with similar arguments for relation \eqref{eq.KPX.3.5-11}, we obtain
\beao
&&\E\left[J_x(T)\,(J_x(T) -1) \right]= \E\left[\sum_{\tau_i\leq T} {\bf 1}_{\left\{ {\bf X}^{(i)}\,e^{-r\,\tau_i}\, \in\; x\,A \right\}}\;\sum_{\tau_j\leq T,\,j\neq i} {\bf 1}_{\left\{{\bf X}^{(j)}\, e^{-r\,\tau_j}\, \in\; x\,A \right\}} \right]\\[2mm]
&&=\E \left[ \E_N\left[\sum_{\tau_i\leq T} {\bf 1}_{\left\{ {\bf X}^{(i)}\,e^{-r\,\tau_i}\, \in\; x\,A \right\}}\;\sum_{\tau_j\leq T,\,j\neq i} {\bf 1}_{\left\{ {\bf X}^{(j)}\,e^{-r\,\tau_j}\, \in\; x\,A \right\}}\;\Big|\;\mathcal{F}_N \right] \right] \\[2mm]
&&=\E \left[ \sum_{\tau_i\leq T} \sum_{\tau_j\leq T,\,j\neq i}  \E_N\left[{\bf 1}_{\left\{ {\bf X}^{(i)}\,e^{-r\,\tau_i} \in x\,A \right\}}\; {\bf 1}_{\left\{ {\bf X}^{(j)}\,e^{-r\,\tau_j}\, \in\; x\,A \right\}}\;\big|\;\mathcal{F}_N\right]  \right] \\[2mm]
&&=\E \left[ \sum_{\tau_i\leq T} \sum_{\tau_j\leq T,\,j\neq i} \PP\left({\bf X}^{(1)}\,e^{-r\,\tau_i}\, \in x\,A\,,\;{\bf X}^{(2)}\, e^{-r\,\tau_j}\, \in x\,A\;\big|\;\mathcal{F}_N \right) \right]  \\[2mm]
&&=\E \left[ \sum_{\tau_i\leq T} \sum_{\tau_j\leq T,\,j\neq i} \PP\left( {\bf X}^{(1)}\,e^{-r\,\tau_i}\, \in\ x\,A\,,\; {\bf X}^{(2)}\,e^{-r\,\tau_j}\, \in x\,A\right) \, \right] \\[2mm]
&&=\E \left[ \int_0^T \int_0^T \PP\left( {\bf X}^{(1)}\,e^{-r\,s}\, \in\ x\,A\,,\; {\bf X}^{(2)}\,e^{-r\,t}\, \in x\,A \right)\,\left(N(ds)\,N(dt) -N(ds)\delta_s(dt) \right) \right] \\[2mm]
&&=\int_0^T \int_0^T \PP\left( {\bf X}^{(1)}\,e^{-r\,s}\, \in\ x\,A\,,\; {\bf X}^{(2)}\,e^{-r\,t}\, \in x\,A \right)\,\alpha^{(2)}(ds,\,dt)\,.
\eeao
From the fact that $(e^{-r\,s},\,e^{-r\,t}) \in [e^{-r\,T},\,1]^2$, the $\{{\bf X}^{(i)}\,,\;i \in \bbn\}$ are $RD_A$, see relation \eqref{eq.KPX.3.5-4}, in combination with Assumption \ref{ass.KPX.3.2} we find from the last relation that
\beam \label{eq.KPX.3.5-15}
\E\left[J_x(T)[J_x(T) -1] \right] \leq C\,K\,\left(\int_0^T \PP\left({\bf X}\,e^{-r\,s} \in x\,A \right) m(ds) \right)^2=C\,K\,\Lambda_x^2\,,
\eeam
is true for any large enough $x>0$. Hence, form relations \eqref{eq.KPX.3.5-14} and \eqref{eq.KPX.3.5-15} we obtain
\beam \label{eq.KPX.3.5-16}
\PP\left( J_x(T) \geq 2 \right) =O(\Lambda_x^2)= o(\Lambda_x)\,,
\eeam
that completes the first step.

For the second step, we take into account Markov inequality to obtain
\beam \label{eq.KPX.3.5-17}
\PP\left( J_x(T) = 1 \right)\leq \PP\left( J_x(T) \geq 1 \right) \leq \E[J_x(T)]= \Lambda_x\,.
\eeam
From the other hand side we find
\beam \label{eq.KPX.3.5-18}
&&\Lambda_x = \sum_{n=0}^{\infty} n\,\PP\left( J_x(T) = n \right)= \PP\left( J_x(T) = 1 \right) +\E\left[ J_x(T)\,{\bf 1}_{\{J_x(T) \geq 2\}} \right] \\[2mm] \notag
&&\leq \PP[J_x(T)=1] + \E\left[ J_x(T)\,\left( J_x(T) -1 \right) \right]=\PP[J_x(T)=1] + o\left(\Lambda_x\right)\,,
\eeam
where at the last step we used relation \eqref{eq.KPX.3.5-15}. Thus, from \eqref{eq.KPX.3.5-17}, \eqref{eq.KPX.3.5-18} we conclude the second step.

Finally, for the last step, we condition with respect to the $\{\tau_i \leq T\}$, and since $e^{-r\,\tau_i} \in \left[e^{-r\,T},\,1 \right]$, for any $i \in \bbn$, from Lemma \ref{lem.KPX.3.2} we obtain
\beao
&&\PP \left( {\bf D}(T) \in x\,A\,,\;J_x(T) =0\;\big|\;\mathcal{F}_N \right) = o\left[\sum_{\tau_i\leq T} \PP\left( {\bf X}^{(i)}\,e^{-r\,\tau_i} \in x\,A \;\big|\;\mathcal{F}_{N} \right)\right] \\[2mm] \notag
&&= o\left[\sum_{\tau_i\leq T} \PP\left( {\bf X}^{(i)}\,e^{-r\,\tau_i} \in x\,A \right)\right]\,.
\eeao
Taking the expectations with respect to $N$, we find
\beam \label{eq.KPX.3.5-19}
&&\PP \left( {\bf D}(T) \in x\,A\,,\;J_x(T) =0 \right) = \E_N\left[ \PP\left( {\bf D}(T) \in x\,A\,,\;J_x(T) =0\;\big|\;\mathcal{F}_N \right) \right] \\[2mm] \notag
&&= \E_N\left[ o\left(\sum_{\tau_i\leq T} \PP\left( {\bf X}^{(i)}\,e^{-r\,\tau_i} \in x\,A \right)\right) \right]=o\left(\E_N\left[ \sum_{\tau_i\leq T} \PP\left( {\bf X}^{(i)}\,e^{-r\,\tau_i} \in x\,A \right)\right] \right) \\[2mm] \notag
&&= o\left(\E_N\left[ \int_{0}^T \PP\left( {\bf X}\,e^{-r\,s} \in x\,A \right)\,N(ds) \right] \right)=o\left(\Lambda_x\right)\,.
\eeam
Relation \eqref{eq.KPX.3.5-19} completes the step 3.

Hence, by \eqref{eq.KPX.3.5-13}, together with the steps $(1) - (3)$ we get \eqref{eq.KPX.3.5-12}.
~\halmos

\noindent{\bf Proof of Corollary \ref{cor.KPX.3.1}}~
From the fact that $MRV \subsetneq \mathcal{S}_A$, from  Theorem \ref{th.KPX.3.1} follows, that it is enough to show
\beao
\int_0^T \PP\left({\bf X}\,e^{-r\,s} \in  x\,A \right)\, m(ds) \sim \mu(A)\,\bV(x)\,\int_{0}^T \,e^{-\alpha\,r\,s}\,m(ds) \,.
\eeao
Indeed, since we have $F \in MRV(\alpha,\,\mu)$, it follows, via the homogeneity of \eqref{eq.KPX.2.7}, the $\mu(A) \in (0,\,\infty)$ and the closure property of $\mathcal{R}_{-\alpha}$ with respect to strong tail equivalence, that $F_A \in \mathcal{R}_{-\alpha}$. So, for the random variable $X_A \stackrel{d}{\sim} F_A$ from \eqref{eq.KPX.2.4} we obtain
\beao
&&\int_0^T \PP\left({\bf X}\,e^{-r\,s} \in  x\,A \right)\, m(ds)  = \int_0^T \PP\left( X_{A}> x\,e^{r\,s} \right)\, m(ds)  \\[2mm]
&&\sim \int_{0}^T \,e^{-\alpha\,r\,s}\,\PP\left( X_{A}> x \right)\, m(ds)  \sim \mu(A)\,\bV(x)\,\int_{0}^T \,e^{-\alpha\,r\,s}\,m(ds)  \,,
\eeao
where at the second step, we applied the dominated convergence theorem, via uniformity of $\mathcal{R}_{-\alpha}$.
~\halmos

Now we need another lemma, that is similar to \cite[Lemma 4.4]{chen:konstantinides:passalidis:2025} and its proof follows by same lines. However, in that paper was used a general stochastic process of logarithmic returns for the investment portfolio, under the framework of renewal driven risk model. For these reasons we give, for sake of completeness the whole proof of the lemma.

\ble \label{lem.KPX.3.4}
Under the conditions of Theorem \ref{th.KPX.3.2}, it holds
\beam \label{eq.KPX.3.13}
 \lim_{M \to \infty} \limsup \dfrac{\sum_{i=M+1}^{\infty} \PP\left({\bf X}^{(i)}\,e^{-r\,\tau_i} \in  x\,A \right)}{ \PP\left({\bf X}^{(1)}\,e^{-r\,\tau_1} \in  x\,A \right)} = 0\,.
\eeam
\ele

\pr~
For some term, let say $i$, in the sum of numerator in \eqref{eq.KPX.3.13}, by \cite[Lemma 1]{yi:chen:su:2011}, we find that for any 
\beao
0< q_1 < J_{F_A}^- \leq J_{F_A}^+ < q_2 < \infty\,,
\eeao 
and for some large enough $x> 0$, there exists a constant $K_1>0$, independent of $\tau_i$, such that it holds
\beao
\PP\left({\bf X}^{(i)}\,e^{-r\,\tau_i} \in  x\,A \right) \leq K_1\,\PP\left({\bf X} \in  x\,A \right)\,\left(\E[e^{-q_1\,r\,\tau_i}]\vee \E[e^{-q_2\,r\,\tau_i} ] \right)\,.
\eeao
Now, choosing the $q_1,\,q_2 $ in such a way that \eqref{eq.KPX.3.3} holds, we
obtain from the last relation and \cite[Theorem 3.3(iv)]{cline:samorodnitsky:1994}, that for some $K_2> 0$ and for some large enough $x>0$ and large enough $i \in \bbn$ we find
\beao
\PP\left({\bf X}^{(i)}\,e^{-r\,\tau_i} \in  x\,A \right) \leq K_1\,K_2\,\PP\left({\bf X}^{(1)}\,e^{-r\,\tau_1} \in  x\,A \right)\,\left(\E[e^{-q_1\,r\,\tau_i}]\vee \E[e^{-q_2\,r\,\tau_i} ] \right)^{1/q}\,.
\eeao
Hence, we obtain
\beao
&&\lim_{M \to \infty} \limsup \dfrac{\sum_{i=M+1}^{\infty} \PP\left({\bf X}^{(i)}\,e^{-r\,\tau_i} \in  x\,A \right)}{ \PP\left({\bf X}^{(1)}\,e^{-r\,\tau_1} \in  x\,A \right)} \\[2mm]
&&\leq K_1\,K_2\,\lim_{M \to \infty} \sum_{i=M+1}^{\infty} \,\left(\E[e^{-q_1\,r\,\tau_i}]\vee \E[e^{-q_2\,r\,\tau_i} ] \right)^{1/q}= 0\,,
\eeao
due to relation \eqref{eq.KPX.3.3}.
~\halmos

\noindent{\bf Proof of Theorem \ref{th.KPX.3.2}}~
For the upper bound in relation \eqref{eq.KPX.3.4} and considering the random variables
\beao
X_A^{(i)}= \sup\{u\;:\;{\bf X}^{(i)} \in u\,A \} = \sup_{{\bf p} \in I_A} {\bf p}^T\,{\bf X}^{(i)}\,,
\eeao
for any $i \in \bbn$,
we have
\beao
&&\PP\left({\bf D}(\infty) \in  x\,A \right)= \PP\left(\sup_{{\bf p} \in I_A} {\bf p}^T\,\left[\sum_{i=1}^{\infty} {\bf X}^{(i)}\,e^{-r\,\tau_i}\right] >  x \right)\\[2mm] \notag
&&\leq \PP\left(\sum_{i=1}^{\infty}  X_{A}^{(i)}\,e^{-r\,\tau_{i}} >  x \right)\sim \sum_{i=1}^{\infty} \PP\left( X_{A}^{(i)}\,e^{-r\,\tau_i} > x \right)\\[2mm] \notag
&&= \sum_{i=1}^{\infty}\PP\left( {\bf X}^{(i)}\,e^{-r\,\tau_{i}} \in x\,A\,,\; \tau_i < \infty \right) = \int_{0}^{\infty} \PP\left({\bf X}\,e^{-r\,s} \in  x\,A \right)\,m(ds)\\[2mm] \,,
\eeao
where at the second step we used the monotonicity of the event sequence
\beao
\left\{ \left(\sup_{{\bf p} \in I_A} {\bf p}^{\top}\,\sum_{i=1}^n {\bf X}^{(i)}\,e^{-r\,\tau_i} >x  \right)\,, n \in \bbn \right\}\,,
\eeao
and $\left\{ \left(\sum_{i=1}^n  X_A^{(i)}\,e^{-r\,\tau_i} >x  \right)\,, n \in \bbn \right\}$  in combination with \cite[Prop. 2.4]{konstantinides:passalidis:2024g} and at the third step we used \cite[Theorem 2]{yi:chen:su:2011}.

Now we deal with the lower bound of relation \eqref{eq.KPX.3.4}. Let $M \in \bbn$. Then we obtain
\beam \label{eq.KPX.3.14}
&&\PP\left({\bf D}(\infty) \in  x\,A \right) \geq \PP\left(\sum_{i=1}^{M} {\bf X}^{(i)}\,e^{-r\,\tau_i} \in  x\,A \right)\\[2mm] \notag
&&\geq \sum_{i=1}^{M} \PP\left( {\bf X}^{(i)}\,e^{-r\,\tau_i}  \in  x\,A \right)- \sum_{1\leq i < j \leq M} \PP\left( {\bf X}^{(i)}\,e^{-r\,\tau_i}  \in  x\,A\,,\;{\bf X}^{(j)}\,e^{-r\,\tau_j}  \in  x\,A \right)\,,
\\[2mm] \notag
&&\sim \sum_{i=1}^{M} \PP\left( {\bf X}^{(i)}\,e^{-r\,\tau_i}  \in  x\,A \right)\,,
\eeam
where at the second step we employ Bonferroni inequality, since the set $x\,A$ is increasing and the summands non-negative, while at the third step we used the fact that the structure of $\{ X_A^{(i)}\,e^{-r\,\tau_i}\,,\;i \in \bbn\}$ is $QAI$  due to \cite[Lemma 3.1]{chen:yuen:2009},( or due to \cite[Theorem 2.2]{li:2013}). From relations \eqref{eq.KPX.3.14} and \eqref{eq.KPX.3.13}, for any $\delta>0$, we can find some large enough $M \in \bbn$, such that it holds
\beam \label{eq.KPX.3.15}
&&\PP\left({\bf D}(\infty) \in  x\,A \right) \gtrsim \left[\sum_{i=1}^{\infty} - \sum_{i=M+1}^{\infty} \right] \PP\left( {\bf X}^{(i)}\,e^{-r\,\tau_i} \in  x\,A \right)\\[2mm] \notag
&&\geq (1-\delta)\,\sum_{i=1}^{\infty} \PP\left( {\bf X}^{(i)}\,e^{-r\,\tau_i}  \in  x\,A \right) = (1-\delta)\,\int_{0}^{\infty} \PP\left({\bf X}\,e^{-r\,s} \in  x\,A \right)\,m(ds) \,.
\eeam
From \eqref{eq.KPX.3.15} and the arbitrariness in the choice of $\delta>0$, we reach the desired lower bound, by letting $\delta \downarrow 0$.
~\halmos

\section{Ruin probabilities} \label{sec.KPX.4}

Now we provide an application of Theorems \ref{th.KPX.3.1} and  \ref{th.KPX.3.2}, on ruin probabilities, over finite and infinite-time horizon, respectively, in a multivariate risk model, under weak enough assumptions on the premiums, and eventual Brownian perturbations. As expected, the asymptotic behavior of the ruin probability coincides with the asymptotic behavior of the discounted aggregate claims, while we can realize the asymptotic insensitivity of the the ruin probability with respect to the Brownian perturbations, which is expected when the claims have heavy-tailed distributions, see for one and multidimensional risk model \cite{chen:wang:wang:2013}, \cite{konstantinides:2025}, \cite{li:2017}, \cite{liu:gao:dong:2024}, among others.

We assume that the insurer's surplus process at the moment $t\geq 0$, is given through the relation
\beam \label{eq.KPX.5.1} \notag
&&{\bf U}(t):= \\[2mm] \notag
&&\left(
\begin{array}{c}
U_{1}(t) \\
\vdots \\
U_{d}(t)
\end{array} \right)=x\,\left(
\begin{array}{c}
l_{1} \\
\vdots \\
l_{d}
\end{array} \right) + \left(
\begin{array}{c}
\int_0^t e^{-r\,s}\,c_{1}(s)\,ds \\
\vdots \\
\int_0^t e^{-r\,s}\,c_{d}(s)\,ds
\end{array} \right) + \left(
\begin{array}{c}
\delta_{1} \\
\vdots \\
\delta_{d}
\end{array} \right)\left(
\begin{array}{c}
\int_0^t e^{-r\,s}\,B_{1}(ds) \\
\vdots \\
\int_0^t e^{-r\,s}\,B_{d}(ds)
\end{array}
\right) -{\bf D}(t)\\[2mm]
&&=: x\,{\bf l} + \int_0^{t} e^{-r\,s}\,{\bf c}(s)\,ds + \vec{\delta} \odot \int_0^{t} e^{-r\,s}\,{\bf B}(ds) -  {\bf D}(t) \,,
\eeam
where $x>0$, represents the insurer's initial capital, which is allocated into $d$-lines of business through vector ${\bf l}=(l_1,\,\ldots,\,l_d)$, with $l_1,\,\ldots,\,l_d > 0$ and $l_1+\cdots + l_d =1$. The vector $\{{\bf c}(t) =(c_1(t),\,\ldots,\,c_d(t)\,)\,,\;t \geq 0 \}$ represents the densities of the premiums of the  $d$-lines of business. Hereafter, we assume that the densities of the premiums are non-negative and bounded from above, namely for any $i=1,\,\ldots,\,d$, there exists some $M_i \in (0,\,\infty)$, such that for any $t \geq 0$ it holds
\beam \label{eq.KPX.5.2}
0\leq c_i(t) \leq M_i \,.
\eeam
The assumption of insurer's bounded densities, provides significant flexibility in premium modeling for the risk model \eqref{eq.KPX.5.1}. Indeed, relation \eqref{eq.KPX.5.2} is a sufficient condition, in order to consider that the premiums are arbitrarily dependent with the discounted aggregate claims ${\bf D}(t)$, which is realistic in actuarial practice. For example some large claim usually affects to future premiums. Further, the condition \eqref{eq.KPX.5.2} is not strict, since allows for any $T \in (0,\,\infty)$ and any $i=1,\,\ldots,\,d$ we obtain
\beam \label{eq.KPX.5.3}
\int_0^T e^{-r\,s}\,c_i(s)\,ds \leq \int_0^{\infty} e^{-r\,s}\,c_i(s)\,ds \leq \dfrac {M_i}r \,,
\eeam
which means that the total discounted premiums are bounded from above
by the total discounted premiums of constant premium density $M_i$, for arbitrarily large, fixed $M_i \in (0,\,\infty)$.

The $\vec{\delta}=(\delta_1,\,\ldots,\,\delta_d) \geq {\bf 0}$, is called diffusion coefficient, while the random field $\{{\bf B}(t) =(B_1(t),\,\ldots,\,B_d(t))\,,\;t\geq 0 \}$ is a multivariate standard Brownian motion, that is independent of all the other sources of randomness, and contains arbitrarily correlated components. Thus, $\{{\bf B}(t)\,,\;t\geq 0 \}$ represents one more source of randomness, that can stem either from the premiums or form the claims, or from any small perturbations in the company by the market. Note that ${\bf B}( 0)={\bf 0}$.

The ruin probability in multivariate risk models, can be defined by several ways,
since the entrance sets of the surplus that indicate some insurer's economic rigidity have not unique definition. We follow Assumption \ref{ass.KPX.4.1} in order to describe the properties of such kind of 'ruin' sets. This assumption was introduced in \cite{samorodnitsky:sun:2016}, and is directly related with set family $\mathscr{R}$. Let us remind that if a set $-\bbb$ is increasing, the set $\bbb$ is called decreasing.

\begin{assumption} \label{ass.KPX.4.1}
Let $L$ be a set that is open, decreasing, ${\bf 0}\in \partial L$, $L^c$ is convex, and
for any $x>0$, it holds $x\,L = L$.
\end{assumption}

\bre \label{rem.KPX.4.1}
It is not difficult to verify that if the set $L$ satisfies Assumption \ref{ass.KPX.4.1}, then $A:={\bf l}- L \in \mathscr{R}$. Through this property follows the direct connection of the distribution class of claim vectors with the ruin sets of interest. Two typical ruin sets $L$, that satisfy Assumption \ref{ass.KPX.4.1} are the following:
\beao
L_1=\left\{{\bf z}\;:\;\sum_{i=1}^d z_i <0 \right\}\,, \qquad L_2=\left\{{\bf z}\;:\;z_i <0\,,\;\exists\;i=1,\,\ldots,\,d \right\}\,,
\eeao
related with the sets $A_1$, $A_2$ from \eqref{eq.KPX.2.2} and \eqref{eq.KPX.2.3}, respectively. So, the entrance probability of surplus process into $L_1$, indicates the probability that at some moment the sum of surpluses of $d$-lines of business becomes negative, while the entrance probability of surplus process into $L_2$, indicates the probability that at some moment the surplus for someone of the $d$-lines becomes negative. We can notice, that the typical one-dimensional ruin set $L_3=(-\infty,\,0)$, satisfies Assumption \ref{ass.KPX.4.1} (because is a special case of $L_2$), and is related to $1- L_3=(1,\,\infty) =A \in \mathscr{R}$ of the one-dimensional tail.
\ere

Now, we can define the ruin probabilities with respect to the ruin sets $L$, of Assumption \ref{ass.KPX.4.1}. For a fixed $T \in (0,\,\infty)$, the finite-time ruin probability in risk model \eqref{eq.KPX.5.1} is expressed as
\beam \label{eq.KPX.5.4}
&&\psi_{{\bf l},L}(x;\,T):=\PP({\bf U}(t) \in L\,,\;\exists\; t \in [0,\,T]) \\[2mm] \notag
&& = \PP\left( {\bf D}(t)- \int_0^{t} e^{-r\,s}\,{\bf c}(s)\,ds - \vec{\delta} \odot \int_0^{t} e^{-r\,s}\,{\bf B}(ds) \in x\,A\,,\;\exists\;t\in [0,\,T] \right)\,,
\eeam
and the corresponding infinite-time ruin probability as
\beao
&&\psi_{{\bf l},L}(x;\,\infty):=\PP({\bf U}(t) \in L\,,\;\exists\; t \geq 0) \\[2mm] \notag
&& = \PP\left( {\bf D}(t)- \int_0^{t} e^{-r\,s}\,{\bf c}(s)\,ds - \vec{\delta} \odot \int_0^{t} e^{-r\,s}\,{\bf B}(ds) \in x\,A\,,\;\exists\;t \geq 0 \right)\,,
\eeao

\bco \label{cor.KPX.4.1}
Let $L$ satisfy Assumption \ref{ass.KPX.4.1}, and $A={\bf l}- L $. We suppose in risk model \eqref{eq.KPX.5.1}
that $\{{\bf B}(t)\,,\;t\geq 0 \}$ is independent from all the other sources of randomness. For any fixed $T \in \Lambda \setminus \{\infty \}$ we have the following:
\begin{enumerate}
\item[(i)]
under conditions of Theorem \ref{th.KPX.3.1} it holds
\beam \label{eq.KPX.5.6}
\psi_{{\bf l},L}(x;\,T)\sim \int_0^{T} \PP({\bf X}\,e^{-r\,s} \in x\,A)\,m(ds) \,,
\eeam
\item[(ii)]
 under conditions of Corollary \ref{cor.KPX.3.1} it holds
\beam \label{eq.KPX.5.7}
\psi_{{\bf l},L}(x;\,T)\sim \mu(A)\,\bV(x)\,\int_0^{T} e^{-\alpha\,r\,s}\,m(ds) \,.
\eeam
\end{enumerate}
\eco

\pr~
\begin{enumerate}
\item[(i)]
For any $i=1,\,\ldots,\,d$, we take
\beam \label{eq.KPX.5.8}
M_i(T) := \delta_i \int_0^{T} e^{-r\,s} \,B_i(ds) \,,
\eeam
and
\beam \label{eq.KPX.5.9}
\underline{M}_i(T) := \inf_{0 \leq t \leq T} M_i(t) \leq 0\,, \qquad \overline{M}_i(T) := \sup_{0 \leq t \leq T} M_i(t) \geq 0\,,
\eeam
(having in mind that $B_i(0)=0$, for any $i=1,\,\ldots,\,d$). The random variables from
\eqref{eq.KPX.5.9} are non-defective, namely the have no mass at $-\infty$, or at $\infty$, respectively, and hence the following quantities
\beam \label{eq.KPX.5.10}
\widehat{\underline{M}}(T) :=\bigwedge_{i=1}^d \underline{M}_i(T) \leq 0\,, \qquad \check{\overline{M}}(T) :=\bigvee_{i=1}^d \overline{M}_i(T) \geq 0\,,
\eeam
are also non-defective random variables, independent of ${\bf D}(T)$ and of the premiums. Therefore, by \eqref{eq.KPX.5.4}, and the fact that the set $x\,A$ is increasing and that ${\bf D}(t)$ is non-decreasing function with respect to $t$, we obtain
\beam \label{eq.KPX.5.11}
&&\psi_{{\bf l},L}(x;\,T) \leq \PP\left({\bf D}(T) - \widehat{\underline{M}}(T) \in x\,A \right) \\[2mm] \notag
&& = \int_0^{\infty} \PP\left({\bf D}(T) +y \in x\,A \right)\,\PP\left( - \widehat{\underline{M}}(T) \in dy\right)\\[2mm] \notag
&& \leq \int_0^{\infty} \PP\left({\bf D}(T) \in (x- k_y)\,A \right)\,\PP\left( - \widehat{\underline{M}}(T) \in dy\right)\\[2mm] \notag
&& \sim \int_0^{\infty} \int_0^T \PP\left({\bf X}\,e^{-r\,s} \in (x- k_y)\,A \right)\,m(ds)\,\PP\left( - \widehat{\underline{M}}(T) \in dy\right)\\[2mm] \notag
&& \sim \int_0^T \PP\left({\bf X}\,e^{-r\,s} \in x\,A \right)\,m(ds)\,,
\eeam
where at the third step we used the fact that from \cite[Lemma 4.3(d)]{samorodnitsky:sun:2016}, for any $y \in [0,\,\infty)$, there exists some $0<k_y$, such that for all $x > k_y$ it holds
\beam \label{eq.KPX.5.12}
(x +k_y)\,A \subsetneq x\,A +y \subsetneq (x - k_y)\,A\,,
\eeam
and at the fourth step we applied Theorem \ref{th.KPX.3.1}. In order to see this, can be used the relation
\beao
\PP\left({\bf D}(T) \in (x- k_y)\,A \right) \sim \int_0^T \PP\left({\bf X}\,e^{-r\,s} \in (x- k_y)\,A \right)\,m(ds)\,,
\eeao
for any $k_y$, that corresponds to some $y\in [0,\,\infty)$. Finally at the last step we took into consideration that $F \in \mathcal{S}_A \subsetneq \mathcal{L}_A$.

Now we deal with the lower bound. Following similar path and using for this case the first inclusion of \eqref{eq.KPX.5.12} we find
\beam \label{eq.KPX.5.13}
&&\psi_{{\bf l},L}(x;\,T) \geq \PP\left({\bf D}(T) - \int_0^{T} e^{-r\,s}\,{\bf c}(s)\,ds - \vec{\delta} \odot \int_0^{T} e^{-r\,s}\,{\bf B}(ds) \in x\,A \right)\\[2mm] \notag
&& \geq \PP\left({\bf D}(T) -\check{\overline{M}}(T) \in (x+ k)\,A \right) =\int_0^{\infty} \PP\left({\bf D}(T) \in (x+ k)\,A+y \right)\,\PP\left( \check{\overline{M}}(T) \in dy\right)\\[2mm] \notag
&& \geq \int_0^{\infty} \PP\left({\bf D}(T) \in (x+ k+k_y)\,A \right)\,\PP\left( \check{\overline{M}}(T) \in dy\right) \\[2mm] \notag
&&\sim \int_0^T \PP\left({\bf X}\,e^{-r\,s} \in x\,A \right)\,m(ds)\,,
\eeam
where here at the second step we also used the condition \eqref{eq.KPX.5.2} of bounded premiums, to find the $k >0$, (recall also \eqref{eq.KPX.5.3}). Hence, from \eqref{eq.KPX.5.11} and \eqref{eq.KPX.5.13} we get relation \eqref{eq.KPX.5.6}.
\item[(ii)]
For the proof we follow the line of part (i), with only exception, instead of Theorem \ref{th.KPX.3.1} we use Corollary \ref{cor.KPX.3.1}.~\halmos
\end{enumerate}

Now, we proceed to the infinite time ruin probability.

\bco \label{cor.KPX.4.2}
Let $L$ satisfy Assumption \ref{ass.KPX.4.1}, and $A={\bf l}- L$. We suppose in risk model \eqref{eq.KPX.5.1} that $\{{\bf B}(t)\,,\;t\geq 0 \}$ is independent from all  the other sources of randomness. Then we obtain the following:
\begin{enumerate}
\item[(i)]
under conditions of Theorem \ref{th.KPX.3.2} it holds
\beam \label{eq.KPX.5.14}
\psi_{{\bf l},L}(x;\,\infty)\sim \int_0^{\infty} \PP({\bf X}\,e^{-r\,s} \in x\,A)\,m(ds) \,,
\eeam
\item[(ii)]
under conditions of Corollary \ref{cor.KPX.3.2} it holds
\beam \label{eq.KPX.5.15}
\psi_{{\bf l},L}(x;\,\infty)\sim \mu(A)\,\bV(x)\,\int_0^{\infty} e^{-\alpha\,r\,s}\,m(ds) \,.
\eeam
\end{enumerate}
\eco

\pr~
\begin{enumerate}
\item[(i)]
We follow similar argumentation with that in proof of Corollary \ref{cor.KPX.4.1}(i).
Hence, we present only the points that are different from there.

We define relations \eqref{eq.KPX.5.8} and \eqref{eq.KPX.5.9} for $T=\infty$ (where $0\leq t \leq \infty$, is understood as $t\geq 0$), and following relation \eqref{eq.KPX.5.10} we find that
\beao
\widehat{\underline{M}}(\infty) :=\bigwedge_{i=1}^d \underline{M}_i(\infty) \leq 0\,, \qquad \check{\overline{M}}(\infty) :=\bigvee_{i=1}^d \overline{M}_i(\infty) \geq 0\,,
\eeao
are non-defective random variables, independent of ${\bf D}(\infty)$ and of the premiums. For the upper bound of \eqref{eq.KPX.5.14}, we repeat the steps of \eqref{eq.KPX.5.11} applying now Theorem \ref{th.KPX.3.2}, to obtain
\beao
&&\psi_{{\bf l},L}(x;\,\infty) \leq \PP\left({\bf D}(\infty) - \widehat{\underline{M}}(\infty) \in x\,A \right) \\[2mm] \notag
&& = \int_0^{\infty} \PP\left({\bf D}(\infty) +y \in x\,A \right)\,\PP\left( - \widehat{\underline{M}}(\infty) \in dy\right)\\[2mm] \notag
&& \leq \int_0^{\infty} \PP\left({\bf D}(\infty) \in (x- k_y)\,A \right)\,\PP\left( - \widehat{\underline{M}}(\infty) \in dy\right)\\[2mm] \notag
&& \sim \int_0^{\infty} \PP\left({\bf X}\,e^{-r\,s} \in x\,A \right)\,m(ds)\,.
\eeao

Now we proceed to the lower bound of \eqref{eq.KPX.5.14}, which follows as in case of \eqref{eq.KPX.5.13} but now applying Theorem \ref{th.KPX.3.2}, to obtain
\beao
&&\psi_{{\bf l},L}(x;\,\infty) \geq \PP\left({\bf D}(\infty) - \int_0^{\infty} e^{-r\,s}\,{\bf c}(s)\,ds - \vec{\delta} \odot \int_0^{\infty} e^{-r\,s}\,{\bf B}(ds) \in x\,A \right)\\[2mm] \notag
&& \geq \PP\left({\bf D}(\infty) -\check{\overline{M}}(\infty) \in (x+ k)\,A \right) =\int_0^{\infty} \PP\left({\bf D}(\infty) \in (x+ k)\,A+y \right)\,\PP\left( \check{\overline{M}}(\infty) \in dy\right)\\[2mm] \notag
&& \geq \int_0^{\infty} \PP\left({\bf D}(\infty) \in (x+ k+k_y)\,A \right)\,\PP\left( \check{\overline{M}}(\infty) \in dy\right) \sim \int_0^{\infty} \PP\left({\bf X}\,e^{-r\,s} \in x\,A \right)\,m(ds)\,,
\eeao
which implies the desired lower bound.
\item[(ii)]
For the proof we follow the line of part (i), with only exception, instead of Theorem \ref{th.KPX.3.2} we apply Corollary \ref{cor.KPX.3.2}.~\halmos
\end{enumerate}

\bre \label{rem.KPX.4.2}
As for the discounted aggregate claims, here again if we assume that the conditions of Corollary \ref{cor.KPX.4.1}(i), (ii) (or Corollary \ref{cor.KPX.4.2}(i), (ii)), become stricter with respect to $F \in \mathcal{S}_{\mathscr{R}}$ and the structure of claims is $RD_A$, for any $A \in \mathscr{R}$, (with respect to $F \in (\mathcal{C}\cap \mathcal{P_D})_{\mathscr{R}}$, the structure of claims is $QAI_A$, for any $A \in \mathscr{R}$, and \eqref{eq.KPX.3.3} is satisfied for any $A \in \mathscr{R}$, respectively), then relations \eqref{eq.KPX.5.6}, \eqref{eq.KPX.5.7} (or, \eqref{eq.KPX.5.14}, \eqref{eq.KPX.5.15}, respectively) are true for any ruin sets $L$, that satisfy Assumption \ref{ass.KPX.4.1}. This shows the generality of the approach in \cite{samorodnitsky:sun:2016}, since in modern risk theory were studied concrete forms of $L$, and was considered only in structure of $MRV$. However, some important ruin sets do not satisfy Assumption \ref{ass.KPX.4.1}, see  in \cite[Remark 2.2]{konstantinides:passalidis:2024g}. In these approaches the $MRV$ structure works, as also in approaches that escape from multivariate linear single big jump principle, see in  \cite{chen:cheng:zheng:2025}, \cite{gao:yang:2014}, \cite{konstantinides:passalidis:2025a}, \cite{yang:su:2023}, \cite{yuan:lu:fu:2025}.
\ere


\begin{thebibliography}{99}



\bibitem{albrecher:asmussen:2006}
{\sc Albrecher, H., Asmussen, S.}\ (2006)
Ruin probabilities and aggregate claims distributions for shot noise Cox processes.
{\em  Scand. Actuar. J.} \textbf{2}, 86--110.


\bibitem{asmussen:schmidli:schmidt:1999}
{\sc Asmussen, S., Schmidli, H., Schmidt, V.}\ (1999)
Tail probabilities for non-standard risk and queueing processes with subexponential jumps.
{\em  Adv. Appl. Probab.} \textbf{31}, 442--447.







\bibitem{bernackaite:siaulys:2017}
{\sc Bernackait\.{e}, E., \v{S}iaulys, J.}\  (2017)
The finite-time ruin probability for an inhomogeneous renewal  risk model.
{\em J. Ind. Managm. Optim.} \textbf{13}, no. 1, 207--222.

\bibitem{bingham:goldie:teugels:1987}
{\sc Bingham. N.H., Goldie, C.M., Teugels, J.L.} \ (1987)
{\em Regular Variation}
Cambridge University Press, Cambridge.

\bibitem{buraczewski:damek:mikosch:2016}
{\sc Buraczewski, D., Damek, E., Mikosch, T.}\ (2016)
{\em Stochastic Models with Power-Law Tails}
Springer, New York.









\bibitem{chen:liu:2022}
{\sc Chen, Y., Liu, J.}\  (2022)
An asymptotic study of systemic expected shortfall and marginal expected shortfall.
{\em Insur. Math. Econom.} \textbf{105}, 238--251.

\bibitem{chen:yuen:2009}
{ \sc Chen, Y., Yuen, K.C.}\ (2009)
Sums of pairwise quasi-asymptotic independent random variables with consistent variation.
{ \em Stochastic Models}, \textbf{25}, 76--89.



\bibitem{chen:cheng:zheng:2025}
{\sc Chen, Z., Cheng, D., Zheng, H.}\ (2025)
On the joint tail behavior of randomly weighted sums of dependent random variables with applications to risk theory.
{\em Scand. Actuar. J.},  1-20

\bibitem{chen:konstantinides:passalidis:2025}
{\sc Chen, Z., Konstantinides, D.G., Passalidis, C.D.}\ (2025)
Asymptotics for aggregated interdependet multivariate subexponential claims with general investment returns.
{\em Preprint, arXiv:2507.23713}.

\bibitem{chen:li:cheng:2023}
{\sc Chen, Z., Li, M., Cheng, D.}\ (2023)
Asymptotics for sum-ruin probabilities of a bidimensional risk model with heavy-tailed claims and stochastic returns.
{\em Stochastics}, \textbf{96}, no.2, 947-967.

\bibitem{chen:wang:wang:2013}
{\sc Chen, Y., Wang, L., Wang, Y.}\ (2013)
Uniform asymptotics for the finite-time ruin probabilities of two kinds of nonstandard bidimensional risk modes.
{\em J. Math. Anal. Appl.}, \textbf{401}, no. 1, 114--129.


\bibitem{chen:yuan:2017}
{\sc Chen, Y., Yuan, Z.}\  (2017)
A rivisit to ruin probabilities in the presence of heavy-tailed insurance and financial risks.
{\em Insur. Math. Econom.} \textbf{73}, 75--81.

\bibitem{cheng:2014}
{\sc Cheng, D.}\  (2014)
Randomly weighted sums of dependent random variables with dominated variation.
{\em J. Math. Anal. Appl.} \textbf{420}, no. 3, 1617--1633.





\bibitem{cheng:konstantinides:wang:2024}
{ \sc Cheng, M., Konstantinides, D.G., Wang, D.}\ (2024)
Multivariate regular varying insurance and financial risks in $d$-dimensional risk model.
{ \em J. Appl. Probab.}, \textbf{61}, no. 4, 1319 -- 1342.


\bibitem{cheng:xu:2020}
{\sc Cheng, F., Xu, H.}\  (2020)
The finite-time ruin probability of the nonhomogeneous Poisson risk model with conditionally independent subexponential claims.
{\em Comm. Stat. Theor. Meth.} \textbf{51}, Vol. 12, 4119--4132.


\bibitem{cline:samorodnitsky:1994}
{\sc Cline, D.B.H., Samorodnitsky, G.}\ (1994)
Subexponentiality of the product of independent random variables.
{\em Stoch. Process. Appl.}, \textbf{49}, 75--98.


\bibitem{dirma:nakiliuda:siaulys:2023}
{\sc Dirma, M., Nakliuda, N., {\v{S}}iaulys, J.}\ (2023)
Generalized moments of sums with heavy-tailed random summands.
{\em  Lith. Math. J.}, \textbf{63}, no. 3, 254--271.




\bibitem{foss:richards:2010}
{\sc Foss, S., Richards, A.} \ (2010)
{\em On sums of conditionally independent subexponential random variables.}
{\em Math. Oper. Resear.}, \textbf{35}, 102--119.



\bibitem{gao:yang:2014}
{\sc Gao, Q., Yang, X.}\ (2014)
Asymptotic ruin probabilities in a generalized bidimensional risk model perturbed by diffusion with constant force of interest.
{\em J. Math. Anal. Appl.}, \textbf{419}, no. 2, 1193--1213.


\bibitem{geluk:tang:2009}
{\sc Geluk, J., Tang, Q.}\ (2009)
Asymptotic tail probabilities of sums of dependent subexponential random variables.
{\em J. Theor. Probab.}, \textbf{22}, 871--882.

\bibitem{geng:liu:wang:2023}
{\sc Geng, B., Liu, Z., Wang, S.}\  (2023)
A Kesten-type inequality for randomly  weighted sums of dependent subexponential random variables with application to risk theory.
{\em Lith. Math. J.} \textbf{63}, 81--91.















\bibitem{jiang:gao:wang:2014}
{\sc Jiang, T., Gao, Q., Wang, Y.}\  (2014)
Max-sum equivalence of conditionally dependent random variables.
{\em Stat. Probab. Lett.}, \textbf{84}, 60--66.

\bibitem{jiang:wang:chen:xu:2015}
{\sc Jiang, T., Wang, Y., Chen, Y., Xu, H.}\  (2015)
Uniform asymptotic estimate for finite-time ruin probabilities of a time-dependent bidimensional renewal model.
{\em Insur. Math. Econom.}, \textbf{64}, 45--53.



\bibitem{ko:tang:2008}
{\sc Ko, B.W., Tang Q.H.} \ (2008)
Sums of dependent non-negative random variables with subexponential tails.
{\em J. Appl. Probab.}, \textbf{45}, 85--94.


\bibitem{konstantinides:2025}
{\sc Konstantinides, D.G.} \ (2026)
Infinite-time ruin probability of a multivariate renewal risk model with Brownian perturbations.
{\em Stat. Probab. Lett.}, \textbf{234}, 110700.

\bibitem{konstantinides:li:2016}
{\sc Konstantinides, D.G., Li, J.}\ (2016)
Asymptotic ruin probabilities for a multidimensional renewal risk model with multivariate regularly varying claims.
{\em Insur. Math. Econom.}, \textbf{69}, 38--44.

\bibitem{konstantinides:mikosch:2005}
{\sc Konstantinides, D.G., Mikosch, T.}\ (2005)
Large Deviations and Ruin Probabilities for Solutions to Stochastic Recurrence Equations with Heavy-tailed Innovations.
{\em Ann. Probab.}, \textbf{33}, 1.992--2.035.

\bibitem{konstantinides:liu:passalidis:2025}
{\sc Konstantinides, D.G., Liu, J., Passalidis, C.D.} \ (2026)
Uniform asymptotics for a multidimensional renewal risk model with multivariate subexponential claims.
{\em Scand. Actuar. J.}, p. 1 -- 21.\\
DOI: 10.1080/03461238.2025.2584008.

\bibitem{konstantinides:passalidis:2025a}
{\sc Konstantinides, D.G., Passalidis, C.D.} \ (2025)
Background risk model in presence of heavy tails under dependence.
{\em Nonlin. Anal. Model. Contr.}, \textbf{30}, no. 5, 982 -- 1010.


\bibitem{konstantinides:passalidis:2024g}
{\sc Konstantinides, D.G., Passalidis, C.D.} \ (2024)
Random vectors in the presence of a single big jump.
{\em Preprint, arXiv:2410.10292}.

\bibitem{konstantinides:passalidis:2025h}
{\sc Konstantinides, D.G., Passalidis, C.D.} \ (2025)
Heavy-tailed random vectors: theory and applications.
{\em Preprint, arXiv:2503.12842}.



\bibitem{konstantinides:passalidis:2024j}
{\sc Konstantinides, D.G., Passalidis, C.D.} \ (2025)
Uniform asymptotic estimates for ruin probabilities of a multidimensional risk model with c\'{a}dl\'{a}g returns and multivariate heavy-tailed claims.
{\em  Insur. Math. Econom.}, \text{125}, 103148.



\bibitem{lehmann:1966}
{\sc Lehmann, E.L.}\ (1966)
Some concepts of dependence.
{\em Ann. Math. Stat.}, \textbf{37}, 1137--1153.




\bibitem{leipus:siaulys:konstantinides:2023}
{\sc Leipus, R., \v{S}iaulys, J., Konstantinides, D.G.}\ (2023)
{\em Closure Properties for Heavy-Tailed and Related Distributions: An Overview.}
Springer Nature, Cham Switzerland.




\bibitem{li:2013}
{\sc Li, J.}\ (2013)
On pairwise quasi-asymptotically independent random variables and  their applications.
{\em Stat. Prob. Lett.}, \textbf{83}, 2081--2087.

\bibitem{li:2016}
{\sc Li, J.}\ (2016)
Uniform asymptotics for a multidimensional time-dependent risk model with multivariate regularly varying claims and stochastic return.
{\em Insur. Math. Econom.}, \textbf{71}, 195--204.

\bibitem{li:2017}
{\sc Li, J.}\ (2017)
A note on the finite-time ruin probabilities of a renewal risk model with Brownian perturbation.
{\em Stat. Prob. Lett.}, \textbf{127}, 49--55.









\bibitem{liu:gao:dong:2024}
{\sc Liu, X., Gao, Q., Dong, Z.}\ (2024)
Asymptotics for a diffusion-perturbed risk model with dependence structures,
constant interest force, and a random number of delayed claims.
{\em Stochastic Models}, \textbf{40}, no. 1, 97--122.
















\bibitem{palmowski:pojer:thonhauser:2025}
{\sc Palmowski, Z., Pojer, S., Thonhauser, S.} \ (2025)
Exact asymptotics of ruin probabilities with linear Hawkes arrivals.
{\em Stoch. Proc. Appl.}, \textbf{182}, 104571.

\bibitem{passalidis:2025}
{\sc Passalidis, C.D.} \ (2025)
Multivariate strong subexponential distributions: Properties and Applications.
{\em Preprint, arXiv:2503.22267}.

\bibitem{resnick:2007}
{\sc Resnick, S.}\ (2007)
{\em Heavy-Tail Phenomena. Probabilistic and Statistical Modeling.}
Springer, New York.



\bibitem{ross:1983}
{\sc Ross, S.M.} (1983)
{\em Stochastic Processes.}
Wiley, New York.




\bibitem{samorodnitsky:sun:2016}
{\sc Samorodnitsky, G., Sun, J.}\ (2016)
Multivariate subexponential distributions and their applications.
{\em Extremes}, \textbf{19}, no. 2, 171--196.




















\bibitem{wang:2011}
{\sc Wang, K.}\ (2011)
Randomly weighted sums of dependent subexponential random variables.
{\em Lith. Math. J.}, \textbf{51}, no. 4, 573--586.


\bibitem{wang:cui:wang:ma:2012}
{\sc Wang, Y., Cui, Z., Wang, K., Ma, X.}\ (2012)
Uniform asymptotics of the finite-time ruin probability for all times.
{\em J. Math. Anal. Appl.}, \textbf{390}, 208--223.

\bibitem{wang:wang:gao:2013}
{\sc Wang, K., Wang, Y., Gao, Q.}\ (2013)
Uniform asymptotics for the finite-time ruin probability of a dependent risk model with a constant interest rate.
{\em Methodol. Comput. Appl. Probab.}, \textbf{15}, 109--124.










\bibitem{yang:su:2023}
{\sc Yang, Y., Su, Q.}\ (2023)
Asymptotic behavior of ruin probabilities in a multidimensional risk model with investment and multivariate regularly varying claims.
{\em J. Math. Anal. Appl.}, \textbf{525}, 127319.


\bibitem{yang:wang:2012}
{\sc Yang, Y., Wang, K.}\ (2012)
Uniform asymptotics for the finite-time and infinite-time ruin probabilities in a dependent risk model with constant interest rate and heavy-tailed claims.
{\em Lith. Math. J.}, \textbf{52}, 111--121.



\bibitem{yang:wang:leipus:siaulys:2013}
{\sc Yang, Y., Wang, K., Leipus, R., \v{S}iaulys, J.}\ (2013)
A note on the max-sum equivalence of randomly weighted sums of heavy-tailed random variables.
{\em Nonlin. Anal. Mod. Contr.}, \textbf{18}, no.4, 519--525.



\bibitem{yi:chen:su:2011}
{\sc Yi, L., Chen Y., Su, C.}\ (2011)
Approximation of the tail probability of randomly weighted sums of dependent random variables with dominated variation.
\emph{ J. Math. Anal. Appl.}, {\bf 376}, 365--372.



\bibitem{yuan:lu:fu:2025}
{\sc Yuan, M., Lu, D., Fu, Y.}\ (2025)
Asymptotics for a multidimensional risk model with a random number of delayed claims and multivariate
regularly varying distribution.
{\em Adv. Appl. Probab.}, 1--30.






\end{thebibliography}
\end{document}